\documentclass[12pt]{article}
\usepackage{array}
\usepackage{amssymb,amsmath,amscd}

\newtheorem{df}{Definition}[section]

\newtheorem{co}{Corollary}[section]
\newtheorem{pp}{Proposition}[section]
\newtheorem{theorem}{Theorem}[section]
\newtheorem{rmk}{Remark}[section]

\let\ssection=\section
\renewcommand{\section}{\setcounter{equation}{0}\ssection}

\def\pg{\varphi_{\ve_{\m{G}^*}}}
\def\pb{\varphi_{\ve_{\m{B}^*}}}
\def\a{\alpha}
\def\b{\beta}
\def\g{\gamma}
\def\w{\wedge}
\def\no{\noindent}
\def\ve{\vert}
\def\o{\otimes}
\def\O{\Omega}
\def\ov{\overline}

\def\lm{\lambda}
\def\d{\delta}
\def\e{\epsilon}

\def\h{\mathcal{SL}_{m \ve n}[[x_{ij}]]}

\def\slmn{\m{SL}_{m\ve n}[[y_{ij},z_{ij}]]}
\def\san{s(\m{AN})[[y_{ij},z_{ij}]]}
\def\sumn{\m{SU}_{m\ve n}[[y_{ij},z_{ij}]]}

\def\n{\nabla}

\def\m{\mathcal}
\def\mf{\mathfrak}
\def\C{\mathbb{C}}
\def\R{\mathbb{R}}

\def\s{SL(m \ve n,\C)}

\def\sl{sl(m \ve n, \C)}

\newcommand{\be}{\begin{equation*}}
\newcommand{\ee}{\end{equation*}}
\newcommand{\ben}{\begin{equation}}
\newcommand{\een}{\end{equation}}
\newcommand{\id}{\relax{\rm 1\kern-.27em l}}

\begin{document}


\vspace{3cm}
\begin{center}
{\Large \bf Super Poisson-Lie structure on $SU(m\ve n)$ via $SL(m\ve n,\C)^{\R}$}\\
[50pt]{\small
{\bf F. Pellegrini}\,\footnote{e-mail:pelleg@iml.univ-mrs.fr}
\\ ~~\\Institut de math\'ematiques de Luminy,
 \\ 163, Avenue de Luminy, 13288 Marseille, France}
\end{center}


\vspace{1cm}

\centerline{\textbf{Abstract}}
\vskip1pc
This paper concerns a super Poisson-Lie structure on the real Lie supergroup $SU(m \ve n)$. In fact, it turns out that  the realification of the complex Lie supergroup $SL(m \ve n, \C)$ is a double of $SU(m \ve n)$ i.e. it is endowed with a structure of super Poisson-Lie which brings down on the supergroup $SU(m \ve n)$. We show that the dual Poisson-Lie supergroup of $SU(m\ve n)$ is $s(AN)$. Reciprocaly, $s(AN)$ inherits a super Poisson-Lie structure from the realification of $SL(m\ve n,\C)$ such that its dual Poisson-Lie supergroup is $SU(m\ve n)$.


\section{Introduction}
A Poisson-Lie group $G$ is a Lie group equipped with a Poisson bracket compatible with the group multiplication $m:G \times G \to G$. It turns out that the Poisson-Lie groups naturally come in dual pairs $(G,\tilde{G})$ where the Poisson-Lie bracket on $G$ ($\tilde{G}$) defines the multiplication on $\tilde{G}$ ($G$). Typically, there may exist several Poisson-Lie structures on a given Lie group $G$ and, therefore, there are several dual pairs $(G,\tilde{G})$ where $G$ is fixed and $\tilde{G}$ varies. The classification of all Poisson-Lie structures was succesfully performed for the case where $G$ is a simple compact group \cite{LS}. It turns out that  (modulo a simple Drinfeld twist of the Cartan subalgebra) there is essentialy a unique Poisson-Lie structure whose corresponding dual group $\tilde{G}=AN$ is given by the Iwasawa decomposition $G^{\C}=GAN$ of the complexified group $G^{\C}$.

The canonical Poisson-Lie structure on the simple compact group $G$ \cite{LW} appears in many applications in mathematics \cite{A}, \cite{Ratiu} and in mathematical physics \cite{Klimcik}. The principal motivation of this paper is to construct its superanalogue, i.e. to consider a simple compact supergroup and to find out  a Poisson bracket on it which would be compatible with the supergroup multiplication. It turns out, that this program is less straightforward than it could seem at the first sight. The principal difficulty emerges at the very beginning and it consists in finding the appropriate definition of the concept of the simple compact supergroup. Roughly speaking, a supergroup (or, rather, the algebra of functions on the supergroup) is a graded commutative Hopf superalgebra. It therefore seems that complex (real) supergroups should be seen as complex (real) superalgebras. In particular, the simple compact supergroup should be certainly  a "real" object and hence tempted to identify it to a certain real Hopf superalgebra. Let us now explain, quite amazingly, that this seemingly natural point of view cannot work. Indeed, having a real graded commutative Hopf superalgebra, we can define its real Lie superalgebra $Lie(G)$. If $Lie(G)$ is a true Lie superalgebra (i.e. the odd part is not empty) then the results of Serganova \cite{Serganova} imply: if $Lie(G)$ is simple then either it does not have the odd part (it is just a Lie algebra) or its even part is not compact. For instance, there is no real simple Lie superalgebra whose even part would be $su(2)$.  

In \cite{Pelleg3}(first section), we argued that the no-go result of Serganova does not imply that the simple compact Lie supergroups do not exists but it rather implies that they should not be viewed as real Hopf superalgebras. Speaking more generally, we found very natural to define "real" Lie supergroups as complex graded commutative Hopf superalgebras equipped with a $*$-structure. The situation here is in some sense analogous to that occuring in the theory of quantum groups where e.g. the compact quantum group $SU_{q}(2)$ is not a real Hopf algebra but it is rather a complex Hopf algebra equipped with a $*$-structure (It turns out that the $*$-real points do not form a real Hopf algebra). Of course, the saying that the real supergroup is compact means that the $*$-structure must have some additional properties explaining the term compact. To this issue is devoted the first part of our thesis \cite{Pelleg3} whose results we explain in section 2.

In \cite{Pelleg3}, we have proposed the axioms which the (compact) $*$-structure on complex supergroup should fulfill and we have proved that every complex supergroup from the series $SL(m \ve n, \C)$ and $OSp(2r,s)$ possesses a unique compact $*$-structure. Then, in \cite{Pelleg2} we have performed a test of plausibility of our definition of the compact simple supergroup. Indeed, we have shown that the realification of the complex supergroup $SL(m \ve n,\C)$  admits a superanalogue of the global Iwasawa decomposition $SL(m \ve n, \C)=SU(m \ve n)s(AN)$ where $SU(m \ve n)$ is our compact real form of $SL(m \ve n,\C)$ and $s(AN)$ is a real supergroup appropriately defined using the positive superroots. Finally, in this paper, we corroborate our definition by showing that our compact simple supergroup $SU(m \ve n)$ can be naturaly equipped with a Poisson-Lie structure coming from the structure of $s(AN)$.

The article is organised as follows: in section 2, we recall the definition of a real Lie supergroup as a complex commutative Hopf superalgebra equipped with a $*$-structure. We also give in this section the definition of real Poisson-Lie supergroup and real Baxter-Lie superalgebra. The crucial point here is that all these structures appear in two versions: the graded and the standard. In section 3, we remind the definition of $SL(m\ve n,\C)^{\R}$, $SU(m\ve n)$ and $s(AN)$ (cf. \cite{Pelleg2} and \cite{Pelleg3}). We establish that $SL(m\ve n)^{\R}$ is the real Drinfeld double of the real supergroup $SU(m \ve n)$ and $s(AN)$, we equip it with a Poisson-Lie structure based on the Yang-Baxter operator $R$ (cf. \cite{STS} adapted to the non super case). Then we show that the maps describing the "embedding" of the real supergroup $SU(m \ve n)$ and $s(AN)$ in its double define also a Poisson ideal with respect to the Poisson-Lie bracket on the double. Factorizing the algebra of "functions" (or, rather, of formal power series) on the double by these Poisson ideals then gives the seeken super Lu-Weinstein Poisson-Lie bracket on the compact supergroup $SU(m \ve n)$ and the supergroup $s(AN)$ such that the dual Poisson-Lie supergroup of $SU(m\ve n)$ is $s(AN)$ and vice-versa. Moreover, we close the article by three supplementaries annexes which give some technical results used in the core of the article.


\section{Real Poisson-Lie supergroup}

In our thesis (cf. \cite{Pelleg3}) we have defined standard (graded) real supergroup. We give here these definitions, theirs infinitesimal counterparts, the definitions of standard (graded) real Baxter-Lie superalgebra as well as the notions of standard (graded) real Poisson-Lie supergoups.\\
\\
First, we recall the notion of supercommutative complex Hopf superalgebra. It is a sextuplets $(H,\mu,1,\Delta,\e,S)$ such that $H=H_{0}+H_{1}$ is a superlinear space, $\mu:H\o H \to H$ such that $\; \mu(f\o g)=(-1)^{\ve f \ve \ve g \ve} \mu(g\o f)$ (supercommutative, for short we note in the sequel $\mu(f\o g)=fg$), $f(gh)=(fg)h$ (associativity), $1f=f1=f$ (unity), the coproduct $\Delta: H \to H \otimes H$ fullfils $(\Delta \otimes) 1 \circ \Delta=(1\otimes \Delta) \circ \Delta$ (coassociativity), the counity $\e:H \to \C$ satisfies $f'\e(f'')=\e(f')f''=f$, the antipode $S:H\to H$ is such that $f'S(f'')=S(f')f''=\e(f)$, $\Delta(fg)=\Delta(f)\Delta(g)$ and  $\e(fg)=\e(f)\e(g)$ with $\forall f,g,h \in H$. Here we have used the sweedler notation $\Delta(f)=f'\o f''$ with $f,f',f'' \in H$. A complex supergroup is a supercommutative complex Hopf superalgebra.\\
\\
Let $H$ be a complex supergroup. A linear map $\d:H \to \C$ is called an $\e$-derivation if it satisfies the following property
\ben
\d(fg)=\e(f)\d(g)+\d(f)\e(g),
\label{epsilonder}
\een
for all $f,g\in H$. We note $\mf{h}$ the set of all $\e$-derivations on $H$. $\mf{h}$ is a complex superlinear space with the following gradation: $\d \in \mf{h}_{0}$ if $\d$ vanishes on $H_{1}$, while $\d \in \mf{h}_{1}$ if $\d$ vanishes on $H_{0}$. This linear space is  endowed with the following superbracket
\be
[\d_{1},\d_{2}](f)=\d_{1}(f')\d_{2}(f'')-(-1)^{\ve \d_{1} \ve \ve \d_{2} \ve}\d_{2}(f')\d_{1}(f''),\quad f \in H.
\ee 
This superbracket is well defined because the superbracket of two $\e$-derivations is again an $\e$-derivation.
It is easy to prove that this superbracket is super-antisymetric i.e.
\be
[\d_{1},\d_{2}]=-(-1)^{\ve \d_{1} \ve \ve \d_{2} \ve}[\d_{2},\d_{1}],
\ee
as well as fullfils the super-Jacobi identity
\be
(-1)^{\ve \d_{1} \ve \ve \d_{3} \ve }[\d_{1},[\d_{2},\d_{3}]]+(-1)^{\ve \d_{3} \ve \ve \d_{2} \ve }[\d_{3},[\d_{1},\d_{2}]]+(-1)^{\ve \d_{2} \ve \ve \d_{1} \ve }[\d_{2},[\d_{3},\d_{1}]]=0.
\ee
Thus ($\mf{h}$,$[.,.]$) is a complex Lie superalgebra, this is the complex Lie superalgebra of the complex supergroup $H$.


\begin{rmk}
The definition of Lie superalgebra of Lie supergroup via the $\e$-derivation is well know in the non-super case, as it is shown in \cite{Waterhouse}.
\end{rmk}

\noindent In view to define real supergroup, we have to introduce the following two kinds of $*$-structures.


\begin{df}
Let $H$ be a Hopf superalgebra. A \underline{standard} $*$-structure of $H$ is an even map $*:H \to H$ such that
\ben
(\Delta x)^{* \o *} = \Delta(x^{*}), \label{p1}
\een
\ben
\e(x^{*}) = \overline{\e(x)},  \label{p2}
\een
\ben
(\lm x + \mu y)^{*} = \bar{\lm} x^{*} + \bar{\mu}y^{*}, \label{p3}
\een
\ben
(x y)^{*} = x^{*} y^{*} , \label{p4}
\een
\ben
S(x^{*})=(S(x))^{*}, \label{p5}
\een
\ben
(x^{*})^{*} =  x, \label{p6}
\een
with $x , y \in H$ and $\lm , \mu \in \C$.\\
If the last property is replaced by the following
\ben
(x^{\bigstar})^{\bigstar} = (-1)^{\ve x \ve} x, \label{p7}
\een
then $\bigstar$ is a \underline{graded} $\bigstar$-structure.
\label{star}
\end{df}



\begin{rmk}
In substance, the graded real structure are already introduced in the papers \cite{GKP}, \cite{SNR1} and \cite{SNR2}. 
\end{rmk}



\begin{df}
1) A \underline{standard real} supergroup is a  supercommutative complex Hopf superalgebra endowed with a \underline{standard $*$-structure}.\\
\\
2) A \underline{graded real} supergroup is a supercommutative complex Hopf superalgebra endowed with a \underline{graded $\bigstar$-structure}. 
\label{strsgg}
\end{df}


\vskip1pc
\noindent 
Now we turn to the infinitesimal version of standard and graded real supergroup.


\begin{df}
1) A \underline{standard real} Lie superalgebra is a complex Lie superalgebra endowed with a \underline{standard real structure}.\\
\\
2) A \underline{graded real} Lie superalgebra is a complex Lie superalgebra endowed with a \underline{graded real structure}.  
\label{strsg}
\end{df}


\noindent  
Thus we have to define the notion of standard and graded real structure.


\begin{df}
Let $\mf{h}$ be a complex Lie superalgebra.\\
1) A \underline{standard real structure} $\phi$ is an even map $\phi:\mf{h} \rightarrow \mf{h}$ such that
\be
\phi(\lambda x + \mu y)= \bar{\lambda} \phi(x) + \bar{\mu} \phi(y),
\ee
\be
\phi^{2}(x)=x,
\ee
\be
\phi([x,y])=[\phi(x),\phi(y)].
\ee
with $\lambda, \mu \in \mathbb{C},\, x,y \in \mf{h}$.\\
\\
2)  A \underline{graded real structure} $\phi$ is an even map $\phi:\mf{h} \rightarrow \mf{h}$ such that
\be
\phi(\lambda x + \mu y)= \bar{\lambda} \phi(x) + \bar{\mu} \phi(y),
\ee
\be
\phi^{2}(x)=(-1)^{\vert x \vert}x,
\ee
\be
\phi([x,y])=[\phi(x),\phi(y)].
\ee 
with $\lambda, \mu \in \mathbb{C},\, x,y \in \mf{h}$.
\end{df}


\begin{rmk}
We frequently note a graded (standard) real Lie superalgebra $\mf{h}$ by the couple $(\mf{h}, \phi)$ where $\phi$ is the graded (standard) real structure of $\mf{h}$. 
\end{rmk}

\begin{rmk}
The automorphisms $\phi$, which fullfil the properties $1)$ or $2)$, have been introduced by Serganova \cite{Serganova}. However, she interprets only the automorphisms of the first kind as real form \cite{Serganova}. In our thesis \cite{Pelleg3} we explain why the second morphisms generate also real form of a complex Lie superalgebra, for this it is crucial to work at the functorial level.
\end{rmk}

\noindent Now, we prove a theorem which associate to a standard (graded) $*$-structure on a supergroup a standard (graded) structure on the Lie superalgebra.


\begin{theorem}
Let ($H$,$*$) be a standard or graded real supergroup. Then the Lie superalgebra $\mf{h}$ is equipped with a standard or graded real structure $\phi$ given by the following formula  
\ben
\phi(\d)(f)=\ov{\d(f^{*})},
\label{phi}
\een
with $\d \in \mf{h}$ and $f \in H$.
\label{str}
\end{theorem}

\noindent
{\bf Proof:} \\
\\
We begin by the proof that $\phi(\d)$ is an $\e$-derivation. Effectively, we have
\be
\phi(\d)(fg)=\ov{\d(f^*g^*)}=\ov{\d(f^*)}\; \ov{\e(g^*)}+\ov{\e(f^*)}\; \ov{\d(g^*)}=\phi(\d)(f)\e(g)+\e(f)\phi(\d)(g).
\ee

\noindent 1) Antilinearity:

\be
\phi(\lm \d_{1}+ \mu \d_{2})(f)=\ov{\lm \d_{1}(f^*)}+\ov{\mu \d_{2}(f^*)}=\bar{\lm}\phi(\d_{1})(f)+\bar{\mu}\phi(\d_{2})(f),
\ee

\no with $\lm, \mu \in \C,\; \d_{1}, \d_{2} \in \mf{h}$ and $f,g \in H$.\\
\\
2) Lie superalgebra morphism:

\begin{eqnarray*}
\phi([ \d_{1} , \d_{2} ])(f) & = &\ov{[ \d_{1} , \d_{2} ](f^*)}\\
&=& \ov{\d_{1}(f'^*)}\; \ov{\d_{2}(f''^*)}-(-1)^{\ve \d_{1} \ve \ve \d_{2} \ve }\ov{\d_{2}(f'^*)}\; \ov{\d_{1}(f''^*)}\\
& = & [ \phi(\d_{1}) , \phi(\d_{2}) ](f), 
\end{eqnarray*}

\no for $\d_{1},\d_{2}\in \mf{h}$ and $f \in H$.\\
\\
3)$\phi^{2}$:\\
\\
Let $*$ be a standard $*$-structure. We have
\be
\phi(\phi(\d))(f)=\ov{\phi(\d)(f^*)}=\ov{\ov{\d((f^*)^*)}}=\d(f)
\ee
with $f \in H$.\\
\\
On the other hand, let $\bigstar$ be a graded $\bigstar$-structure, we have 
\be
\phi(\phi(\d))(f)=\ov{\phi(\d)(f^*)}=\ov{\ov{\d((f^*)^*)}}=(-1)^{\ve f \ve}\d(f)=(-1)^{\ve \d \ve} \d(f)
\ee
$\d \in \mf{h}, \quad f \in H.$ 
\\
\\
This ends the proof.\hfill $\blacksquare$
\vskip2pc
\no In view to illustrate these definitions, we study the complex supergroup $SL(m\ve n,\C)$. First, we define $\C[[x_{ij}]]$. It is the superalgebra of formal series generated by $x_{ij}$ for $i,j=1...m+n$, it is a supercommutative Hopf superalgebra for the coproduct, counity and antipode 
\be
\Delta(x_{ij})=1 \o x_{ij}  + x_{ij} \o 1 + \sum_{k=1}^{m+n} x_{ik} \o x_{kj},\quad \e(x_{ij})=0, 
\ee  
\be
S(x_{ij})=-\delta_{ij}+(1 + X)_{ij}^{-1},
\ee
where $X$ is the matrice defined by $(X)_{ij}=x_{ij}$ and $(1 + X)^{-1}$
is the inverse of the sum of the unity matrice $1$ with the matrice $X$. These three maps are defined on all the elements of $\C[[x_{ij}]]$ because they are superalgebra morphisms. The gradation of the generators $x_{ij}$ is $\ve x_{ij} \ve=\ve i \ve + \ve j \ve$ where $\ve i \ve=0$, $\ve j\ve=1$ for respectively $i=1...m, j=m+1...m+n$. Furthermore, the generators fullfil the following equalities $x_{ij}x_{kl}=(-1)^{(\ve i \ve + \ve j \ve)(\ve k \ve +\ve l \ve)}x_{kl}x_{ij}$ so that $\C[[x_{ij}]]$ is supercommutative. The counity and the coproduct seem surely more familiar to the reader on the generators   $u_{ij}=\delta_{ij}+x_{ij}$ i.e.
\be
\Delta(u_{ij})=\sum_{k=1}^{m+n}u_{ik} \o u_{kj}, \quad \e(u_{ij})=\delta_{ij}.
\ee
\\
The supergroup $\s$ \; is the supercommutative Hopf superalgebra $\C[[x_{ij}]]$ factorised by the ideal of Hopf superalgebra generated by the relation $sdet(1+X)-1=0$ with the superdeterminant defined by
\be
sdet\left(
\begin{array}{cc}
       A & B   \\
       C & D
\end{array}
\right)=\frac{det(A-BD^{-1}C)}{det(D)}
\ee
where $A$, $B$, $C$, $D$ are respectively $m\times m$-matrices, $m\times n$-matrices, $n\times m$-matrices and $n\times n$-matrices. We note the complex Hopf superalgebra of $SL(m\ve n,\C)$ by  $\h$.\\ 
\\ 
Now, we turn to the Lie superalgebra of $\h$, we note it by $sl(m\ve n,\C)$. By definition $sl(m\ve n,\C)$ is the set of $\e$-derivations on $\C[[x_{ij}]]$ which vanish on the ideal generated by $sdet(1+X)-1=0$. Thus, we begin by determining the space of $\e$-derivations on $\C[[x_{ij}]]$. In order to know an $\e$-derivation on $\C[[x_{ij}]]$, it is sufficient to evaluate it on the generators $x_{ij}$. We note $\d(x_{ij})=M_{ij}$, thereby $\d$ defines a supermatrice\footnote[2]{The set of complex  $(m+n) \times (m+n)$ supermatrices, noted $M(m \ve n,\C)$, forms a superlinear space which elements are such that
$$M=\left(
\begin{array}{cc}
       P & Q   \\
       R & T
\end{array}
\right)$$
where $P,T$ are respectively $m\times m$, $n \times n$ complex matrices and $Q,R$ are respectively $m\times n$, $n \times m$ complex matrices.
The even complex matrices are such that $Q=R=0$ and odd complex supermatrices satisfie $P=T=0$. The superbracket of $M(m \ve n,\C)$ is $[M,N]=MN-(-1)^{\ve M \ve \ve N \ve}NM$.} $M$ which have the same parity of $\d$. Reciprocaly any supermatrice $M$ define an $\e$-derivation by the following formula $\d_{M}(x_{ij})=M_{ij}$. Moreover, we have the following equality
\be
[\d_{M},\d_{N}](x_{ij})=\d_{[M,N]}(x_{ij}).
\ee
Hence, the complex Lie superalgebra of $\C[[x_{ij}]]$ is isomorph to the Lie superalgebra of supermatrices $M(m\ve n,\C)$. Moreover we have 
\be
\d_{M}(sdet(1+X))=-trT+trP=StrM.
\ee
Thereby $sl(m\ve n,\C)=\{M\in M(m\ve n,\C)/ Str(M)=0 \}$.\\
\\
The couple $(\h,\bigstar)$ is a graded real supergroup with
\be
x^{\bigstar}_{ij}=(-1)^{(\ve i \ve + \ve j \ve)\ve j \ve}S(x_{ji}).
\ee
The proof that this $\bigstar$-structure is graded is done in \cite{Pelleg2} and \cite{Pelleg3}. 
This $\bigstar$-structure gives on the Lie superalgebra $sl(m\ve n,\C)$ the following graded real structure
\be
\phi(\d_{M})(x_{ij})=\d_{\phi(M)}(x_{ij})=\overline{\d_{M}(x^{\bigstar}_{ij})}=\d_{(-1)^{\ve M \ve}\bar{M}^{st}}(x_{ij}).
\ee
Thus $\phi(M)=(-1)^{\ve M \ve}\bar{M}^{st}$ where the supertranspose of a supermatrice is defined by the following equality
\be
\left(
\begin{array}{cc}
P & Q \\
R & T
\end{array}
\right)^{st}=
\left(
\begin{array}{cc}
P^t & R^t \\
-Q^t & T^t
\end{array}
\right),
\ee
with $t$ the usual transposition of matrices.
We remark that the real Lie algebra of fixed points of $\phi$ is $su(m)\oplus su(n) \oplus u(1)$ which is the direct  sum of three compact real forms. Thereby, we say that $(\h,\bigstar)$ is a compact graded real supergroup and we call it $SU(m\ve n)$. 
\vskip5pc
\no We have previously defined real supergroup like supercommutative Hopf superalgebra with a $*$-structure. Thus, if we endow a real supergroup with a supplementary structure (for instance a Poisson superbracket), then it is necessary to impose  that the $*$-structure and the extra structure satisfy some compatibility relation. Here is an example with real Poisson-Lie supergroup.


\begin{df}
Let ($H$,*) a standard (or graded) real supergroup. ($H$,\{.,.\},*) is a  standard (or graded) real Poisson-Lie supergroup if it exists a bilinear map $H\times H \to H$ noted $(f,g)\to\{f,g\}$ such that for all $f,g,h\in H$\\
\\
i) it fullfils the super-Jacobi identity: 
\be
(-1)^{\ve f \ve \ve h \ve}\{f,\{g,h\}\}+(-1)^{\ve h \ve \ve g \ve}\{h,\{f,g\}\}+(-1)^{\ve g \ve \ve f \ve}\{g,\{h,f\}\}=0,
\ee 
ii) it is  super-antisymetric:
\be
\{f,g\}=-(-1)^{\ve f \ve \ve g \ve}\{ g,f\}, 
\ee
iii) it fullfils the super-Leibniz rule:
\be
\{f,gh\}=\{f,g \}h+(-1)^{\ve f \ve \ve g \ve}g\{f,h\},
\ee
iv) the coproduct $\Delta$ of $H$ is a Poisson morphism:
\be
\Delta\{ f,g\}=\{\Delta(f),\Delta(g)\},
\ee
v) it is compatible with the $*$-structure:
\be
\{f^*,g^*\}=\{ f,g\}^*.
\ee
\label{defspl}
\end{df}

\begin{rmk}
When a map $(f,g)\to \{f,g\}$ fullfils the properties $i)$ to $iii)$, one said that it is a Poisson superbracket, if furthermore it fullfils $iv)$ one said that it is a Poisson-Lie superbracket. The superbracket is defined on $H\o H$ (cf. $\cite{Andru}$) by
\be 
\{f_{1}\o f_{2},g_{1}\o g_{2}\}=(-1)^{\ve g_{1} \ve \ve f_{2} \ve}\{f_{1},g_{1}\}\o f_{2}g_{2}+(-1)^{\ve g_{1} \ve \ve f_{2} \ve} f_{1}g_{1}\o \{f_{2},g_{2}\}.
\ee
It fullfils automaticaly the super-Jacobi identity,  super-Leibniz rule and superantisymetry if $\{.,.\}$ fullfils them.
\end{rmk}

\begin{rmk}
The  Poisson-Lie groups and quantum groups are intimetely linked, in fact the first are obtained as quasi-classic limit of the last. In particular, the compatibility  of the product and the $*$-structure gives precisely, via the quasi-classic limit, our condition $v)$, which exprims the compatibility of the Poisson superbracket with the $*$-structure.
\end{rmk}

\no In the same order of idea, we give the definition of a Poisson-Lie sub-supergroup.


\begin{df} 
Let ($H$,\{.,.\},*) a standard (or graded) real Poisson-Lie supergroup. Then $K$,  is a Poisson-Lie sub-supergroup of ($H$,\{.,.\},*) if \\
\\
i) it exists an ideal of Hopf superalgebra $I$ such that $K=H/I$,\\
ii) $I^*\subset I$,\\
iii) $\{I,H\}\subset I$ .
\label{ssgpl}
\end{df}

\begin{rmk}
An ideal $I$ of Hopf superalgebra $H$ is a subset of $H$ such that $I.H\subset I$, $\e(I)=0$, $\Delta(H)\subset I\o H \dotplus H \o I$ and $S(I)\subset I$. \\
\\
The property $iii)$ in the previous definition means that $I$ is an ideal of  Poisson superalgebra.
\end{rmk}

\no As in the non-super case, a super Poisson-Lie structure equips the dual linear of its Lie superalgebra with a structure of Lie superalgebra. Before proving this result, we introduce the linear dual of the Lie superalgebra of a supergroup. In fact, we have:


\begin{theorem}
Let $H$ be a complex supergroup and $Ker\e$ be the kernel of the counit $\e$ of $H$. Then $\mf{m}=Ker\e/(Ker\e)^2$ is the linear dual of the Lie superalgebra $\mf{h}$ of $H$ i.e $\mf{m}^*=\mf{h}$.
\end{theorem}


\no{\bf Proof:}\\
\\
Let $\d$ be an element of $\mf{h}$ and we note by the same letter its restriction to $Ker\e$. As $\d$ fullfils the property $\d(fg)=\d(f)\e(g)+\e(f)\d(g), \; \forall f,g \in H$, we deduce that $\d$ vanish on $m^{2}_{\e}$. So $\d$ factorise through a linear map  $[\d]:\mf{m} \to \C$ i.e. $[\d] \in \mf{m}^{*}$. Thus, the duality between $\mf{h}$ and $\mf{m}$ is given by the following formula
\ben
\langle [\d], \Omega(f) \rangle=\d(f)
\label{dualite}
\een
where $\Omega(f)$ means the equivalence classes of $f$ in $\mf{m}$ for $f\in Ker\e$.\\
\\
Reciprocaly, let $d \in \mf{m}^{*}$ and $\d$ be a linear map from $H$ to $\C$ such that $\d(f)=d([f-\e(f)])$. Since, we have
\be
fg-\e(fg)=(f-\e(f))\e(g)+(g-\e(g))\e(f)+(f-\e(f))(g-\e(g)), 
\ee
i.e. $fg-\e(fg)$ and $(f-\e(f))\e(g)+(g-\e(g))\e(f)$ define the same equivalence classes, we deduce
\be
\d(fg)=d([fg-\e(fg)])=d([f-\e(f)])\e(g)+\e(f)d([g-\e(g)])=\d(f)\e(g)+\e(f)\d(g).
\ee
These two linear morphisms are clearly inverse of each other. Moreover, $\mf{m}^*$ is equipped of the following superbracket 
\be
[x,y]_{\mf{m}^*}([f])=x([f'-\e(f')])y([f''-\e(f'')])-(-1)^{\ve x \ve \ve y \ve}y([f'-\e(f')])x([f''-\e(f'')]
\ee
for all $f\in Ker\e$ and $x,y \in \mf{m}^*$. And it is easy to observe that for the two previous maps the image of the superbracket is the superbracket of the image, thus $\mf{h}$ and $\mf{m}^*$ are isomophic Lie superalgebras.
\hfill$\blacksquare$
\vskip2pc



\begin{rmk}
This theorem is a generalisation of a well known theorem in the non-super case, as we can find it in \cite{Perrin}. We give it in order to keep this paper self-contained.
\end{rmk}


\no As in the non-super case, a Poisson-Lie structure induces on the linear dual of its Lie superalgebra a structure of Lie superalgebra, this is the following proposition.


\begin{pp}
Let $(H,\bigstar,\{.,.\}_{H})$ be a graded (or standard) real Poisson-Lie supergroup and $(\mf{h},\phi)$ its graded (or standard) real Lie superalgebra.
Then the linear dual $\mf{h}^*=Ker\e_{H}/(Ker\e_{H})^2$ of $\mf{h}$ is canonicaly endowed with a graded (or standard) real Lie superalgebra given by the following formula
\ben
[\O_{H}(f),\O_{H}(g)]_{\mf{h}^*}=\O_{H}(\{f,g\}_{H}),\quad \varphi(\O_{H}(f))=\O_{H}(f^\bigstar)
\een
with $\O_{H}:Ker\e_{H} \to Ker\e_{H}/(Ker\e_{H})^2$ and $f,g\in Ker\e_{H}$.
\label{pp1}
\end{pp}



\no{\bf Proof:}\\
\\
The gradation of $\mf{h}^*$ is the following: $\O_{H}(f)\in \mf{h}^*_{0}$ if $f\in Ker\e_{H}$ and $\ve f \ve=0$ whereas $\O_{H}(f)\in \mf{h}^*_{1}$ if $f\in Ker\e_{H}$ and $\ve f \ve=1$. As the coproduct of $H$ is a Poisson morphism, we deduce that $\e_{H}(\{f,g\}_{H})=0, \; \forall f,g \in H$. Thus, for $f,g\in Ker\e_{H}$ we have $\{f,g\}_{H}\in Ker\e_{H}$, so $\{.,.\}_{H}$ is defined on $Ker\e_{H}$. Moreover, as $\{.,.\}_{H}$ fullfils the super-Leibniz rule, it turns out that $(Ker\e_{H})^2$ is a Poisson ideal of the restriction of $\{.,.\}_{H}$ to $Ker\e_{H}$. Then, it is easy to observe that the following bilinear map $[.,.]_{\mf{h}^*}\mf{h}^*\times \mf{h}^*\to \mf{h}^*$ defined by
\be
[\O_{H}(f),\O_{H}(g)]_{\mf{h}^*}=\{\O_{H}(f),\O_{H}(g)\}_{H}=\O_{H}(\{f,g\}_{H})
\ee 
is super-antisymetric and satisfies the super-Jacobi identity since  $\{.,.\}_{H}$ fullfils them.\\
\\
On the other hand, the map $\varphi:\mf{h}^*\to \mf{h}^*$ such that $\varphi(\O_{H}(f))=\O_{H}(f^\bigstar)$ is well defined i.e. $f^\bigstar \in Ker\e_{H}$ because $\e_{H}(f^\bigstar)=\overline{(\e_{H}(f))}$. It is clearly antilinear and involutive since 
$(\lm f+\mu g)^\bigstar=\bar{\lm}f^{\bigstar}+\bar{\mu}g^{\bigstar}$ and $(f^\bigstar)^\bigstar=(-1)^{\ve f \ve}f$ $\; \forall f,g\in H$ and $\lm , \, \mu \in \C$. Furthermore, from $\{f^\bigstar,g^\bigstar\}_{H}=\{f,g\}^\bigstar_{H}$ we deduce that $\varphi$ is also a morphism of Lie superalgebra. To conclude, $(\mf{h}^*,\varphi)$ is a graded real Lie superalgebra. The standard case follows from similar arguments. \hfill$\blacksquare$
\vskip2pc


\no It exists many ways to equip a Lie group $G$ with a Poisson-Lie bracket.
Among them, there are two which are frequentely used. The first uses a structure of Baxter-Lie algebra on the Lie algebra of the Lie group (it is the approach of the classical $r$-matrix). The second is based  on the injection of the Lie group $G$ in a bigger group $D$, which is already equipped with a structure of Poisson-Lie group. The injection must fullfil the axioms of the non-super version of the definition $\ref{ssgpl}$ i.e. the structure of Poisson-Lie on $G$ comes from the Poisson-Lie group $D$. In the next section, we used the two points of views. First we construct a classical $r$-matrix which endows $SL(m\ve n,\C)^{\R}$ with a Poisson-Lie structure and secondly we show that this latter go down on $SU(m\ve n)$ such that it becomes a Poisson-Lie supergroup.\\
\\
It remains to define the notion of Baxter-Lie superalgebra (in the non-super case see \cite{STS}).
 

\begin{df}
Let $(\mf{h},\phi)$ be a standard (graded) real Lie superalgebra (cf. def. $\ref{strsg}$), $R$ an even linear map on $\mf{h}$. We said that $(\mf{h},\phi,R)$ is a standard (graded) real Baxter-Lie superalgebra if\\
\\
i) $\mf{h}$ is provided with an invariant scalar product noted $(.,.)_{\mf{h}}$ such that
\be
(\phi(x),\phi(y))_{\mf{h}}=\overline{(x,y)}_{\mf{h}},\quad \forall x,y \in \mf{h},
\ee
\\
ii) $R$ is antisymetric:
\be
(R(x),y)_{\mf{h}}=-(x,R(y))_{\mf{h}},\; \forall x,y \in \mf{h},
\ee
satisfies the Baxter-Lie equation:
\be
[R(x),R(y)]=R([R(x),y]+[x,R(y)])-[x,y],\quad \forall x,y \in \mf{h}.
\ee
and fullfils the following relation of compatibility with the standard (graded) real structure $\phi$:
\be
\phi(R(x))=R(\phi(x)),\quad x\in \mf{h}. 
\ee
\end{df}

\begin{rmk}
Thus, in the non-super case, the relations of compatibilities of the real structure $\phi$ with $R$ and the scalar product allow to these last two maps to be defined on the space of fixed points of $\phi$. In other words, the real form associated to $\phi$ is a Baxter-Lie algebra on $\R$.     
\end{rmk}


\section{Super Poisson-Lie structure on $SU(m\ve n)$ via $SL(m\ve n,\C)^{\R}$}

\no We begin with the definition of the graded real supergroups $SL(m\ve n,\C)^{\R}$, $SU(m\ve n)$ and $s(AN)$, where it clearly appears that $SU(m\ve n)$ and $s(AN)$ are sub-supergroups of $SL(m\ve n,\C)^{\C}$ (cf. $\cite{Pelleg2}$ and $\cite{Pelleg3}$). Next, we determine their Lie superalgebra, we endow the Lie superalgebra of $SL(m\ve n,\C)^{\R}$ with a structure of graded real Baxter-Lie superalgebra. Finaly, we establish that $SL(m\ve n,\C)^{\R}$ is the double of $SU(m\ve n)$ and $s(AN)$, these two last supergroups becoming dual graded real Poisson-Lie supergroup. 


\begin{df}
The Hopf superalgebra and $\bigstar$-structure which define the graded real supergroup 
$SL(m\ve n,\C)^{\R}$ are respectively
\be
\slmn=\m{SL}_{m\ve n}[[y_{ij}]]\o \m{SL}_{m\ve n}[[z_{ij}]]
\ee
\ben
y^\bigstar_{ij}=(-1)^{(\ve i \ve + \ve j \ve)\ve j \ve }S(z_{ji}),\quad z^\bigstar_{ij}=(-1)^{(\ve i \ve + \ve j \ve)\ve j \ve }S(y_{ji}).
\label{newstar}
\een
\end{df}


\begin{df}
The Hopf superalgebra which define the graded real supergroup 
$SU(m\ve n)$ is 
\be
\sumn=\slmn/I
\ee
where $I$ is the Hopf ideal generated by the relations
\be
y_{ij}-z_{ij}=0,\quad \forall i,j=1...m+n.
\ee
The ideal $I$ fullfils $I^{\bigstar}\subset I$, where $\bigstar$ is the graded $\bigstar$-structure of $\slmn$ (cf. eq. $(\ref{newstar})$), thus $\sumn$ is endowed with a graded $\bigstar$-structure.
\end{df}


\begin{df}
The Hopf superalgebra which define the graded real supergroup 
$s(AN)$ is 
\be
\slmn/J
\ee
where $J$ is the Hopf ideal generated by the relations 
\be
y_{ji}=z_{ij}=0\; \forall i>j,\quad S(y_{ii})-z_{ii}=0.
\ee
The ideal $J$ fullfils $J^{\bigstar}\subset J$, where $\bigstar$ is the graded $\bigstar$-structure of $\slmn$ (cf. eq. $(\ref{newstar})$), thus $\sumn$ is endowed with a graded $\bigstar$-structure. \end{df}


\begin{rmk}
In the sequel, we work preferentialy with the generators  $u_{ij}=\delta_{ij}+y_{ij}$ and $w_{ij}=\delta_{ij}+z_{ij}$. One remarks that in these variables the coproduct and conunit are 
\be
\Delta(u_{ij})=\sum_{k=1}^{m+n}u_{ik}\o u_{kj},\; \Delta(w_{ij})=\sum_{k=1}^{m+n}w_{ik}\o w_{kj},\; \e(u_{ij})=\d_{ij},\; \e(w_{ij})=\d_{ij}.
\ee
\end{rmk}


\vskip2pc
\no Now, we determine the Lie superalgebra of these supergroups. \\
\\
The Lie superalgebra $\mf{d}$ of $\slmn$ consists in the $\e$-derivations
(cf. $(\ref{epsilonder})$)
$\d_{(A,B)}$ such that
\be
\d_{(A,B)}(y_{ij})=A_{ij},\quad \d_{(A,B)}(z_{ij})=B_{ij},\quad \forall A,B \in \sl.
\ee 
Thus, $\mf{d}$ is isomorph to $\sl \oplus \sl$. Furthermore, the $\bigstar$-structure of $\slmn$ (cf. eq. $(\ref{newstar})$) provides $\mf{d}$ with the following graded real structure (cf. th. $\ref{str}$)
\be
\phi(\d_{(A,B)})(f)=\d_{\phi(A,B)}(f)=\ov{\d_{(A,B)}(f^\bigstar)},\quad \phi(A,B)=(-(-1)^{\ve B \ve}B^{st},-(-1)^{\ve A \ve}A^{st})
\ee
because one have
\be
\ov{\d_{(A,B)}(y^\bigstar_{ij})}=-(-1)^{(\ve i \ve + \ve j \ve)\ve j \ve}\bar{B}_{ji},\quad \ov{\d_{(A,B)}(z^\bigstar_{ij})}=-(-1)^{(\ve i \ve + \ve j \ve)\ve j \ve}\bar{A}_{ji}.
\ee
\vskip1pc
\no 
The Lie superalgebra $\mf{g}$ of $\sumn$ is composed of the $\e$-derivations of $\slmn$ which vanish on the Hopf ideal $I$ i.e. $\mf{g}=\{\d_{(A,A)},\; A \in \sl \}$. Moreover, $\mf{g}$ is a graded real Lie superalgebra because the graded real structure $\phi$  of $\mf{d}$ fullfils $\phi(\mf{g})\subset \mf{g}$ and we note this graded real Lie superalgebra by $(\mf{g},\phi_{G})$ where $\phi_{G}=\phi_{\ve_{\mf{g}}}$.\\
\\
The Lie superalgebra $\mf{b}$ of $\san$ is the set of $\e$-derivations of \\ $\slmn$ which vanish on the Hopf ideal $J$. It is easy to prove that $\mf{b}=\{\d_{(A,B)}/ A\in (\sl_{-} \oplus \sl_{0}),\; B\in (\sl_{0} \oplus \sl_{+}),\; A_{0}+B_{0}=0 \}$ (for the notations cf. Annexe $1$ on $\sl$). Furthermore, $\mf{b}$ is  a graded real Lie superalgebra because the graded real structure $\phi$ fullfils $\phi(\mf{b})\subset \mf{b}$ and we note this graded real Lie superalgebra by $(\mf{b},\phi_{B})$ where $\phi_{B}=\phi_{\ve_{\mf{b}}}$.\\
\\
Now, we show that $\mf{d}$ is graded real Baxter-Lie superalgebra. Indeed, $\mf{d}$ is endowed with the following supersymetric invariant scalar product 
\ben
((A,B),(C,D))_{\mf{d}}=(A,C)_{sl}-(B,D)_{sl},
\label{prdscs}
\een 
where $(A,C)_{sl}=-\frac{i}{2}Str(AC)$ (see Annexe $1$).
This scalar product fullfils the relation of compatibility with the graded real structure $\phi$ i.e. one have 
\be
(\phi(A,B),\phi(C,D))_{\mf{d}}=\ov{((A,B),(C,D))}_{\mf{d}}.
\ee
$\mf{d}$ is also provided with the linear operator $R_{\mf{d}}$ defined by
\be
R_{\mf{d}}=P_{\mf{b}}-P_{\mf{g}}.
\ee
This operator is antisymetric, satisfies the Yang-Baxter equation and the compatibility relation with $\phi$.


\begin{rmk}
One observe that $\mf{g}$ and $\mf{b}$ are isotropic for the scalar product $(\ref{prdscs})$, $\mf{d}=\mf{g}\dotplus \mf{b}$ and $\mf{g}$, $\mf{b}$ are Lie sub-superalgebra of $\mf{d}$, nevertheless $[\mf{g},\mf{b}]\neq 0$. These properties imply that $(\mf{d},\mf{g},\mf{b})$ is a Manin supertriple as it is defined in $\cite{Andru}$.   
\end{rmk}
 
\no For a classical $r$-matrix $R$ on a Lie algebra of a Lie group $G$, it is associated in $\cite{STS}$ a Poisson-Lie bracket on $G$. The super case have been traited in $\cite{Andru}$, thus the superbracket $(\ref{crochetPL})$ in the following theorem is inspired from the one defined in $\cite{Andru}$.\\
\\ 
It is time to write the mains theorems of this article. 


\begin{theorem}
\no Let $SL(m\ve n,\C)^{\R}$ be the graded real supergroup equipped with the following superbracket 

\ben
\{ f,g \}=\frac{1}{2}\sum_{a=1}^{m+n}(-1)^{\ve f \ve \ve a \ve}
(R_{\mf{d}}(\hat{h}_{a}),\hat{h}_{b})_{\mf{d}}
[\n^L_{h_{a}}f \; \n^L_{h_{b}}g-\n^R_{h_{a}}f \; \n^R_{h_{b}}g], 
\label{crochetPL}
\een

\no $\forall f,g \in \slmn$, where $h_{a}$ and $\hat{h}_{b}$ are the dual basis of $\mf{d}$ i.e. $(h_{a},\hat{h}_{b})_{\mf{d}}=\delta_{ab}$. Then we have the following properties:\\
\\
i) this superbracket is a Poisson-Lie superbracket on $(\slmn,\bigstar)$,\\
\\
ii) $I$ and $J$ are Poisson ideals,\\
\\
iii) $\{ f^\bigstar,g^\bigstar \}=\{ f,g\}^\bigstar$ where $y^\bigstar_{ij}=(-1)^{(\ve i \ve + \ve j \ve) \ve j \ve}S(z_{ji}),\; z^\bigstar_{ij}=(-1)^{(\ve i \ve + \ve j \ve) \ve j \ve}S(y_{ji})$. 
\label{th1}
\end{theorem}


\begin{co}
$\sumn$ and $\san$ are graded real Poisson-Lie sub-supergroups of $\slmn$  for the previous superbracket $(\ref{crochetPL})$.
\label{co1} 
\end{co}


\begin{theorem}
Let $(SU(m\ve n),\bigstar)$ and $(s(AN),\bigstar)$ be the graded real Poisson-Lie supergroup defined in the previous corollary, $(\mf{g},\phi_{G})$ and $(\mf{b},\phi_{B})$ its graded real Lie superalgebras. We note by $\mf{g}^*$ and $\mf{b}^*$ the linear supervector spaces of respectively $\mf{g}$ and $\mf{b}$. 
Then the Poisson-Lie superbracket and the graded $\bigstar$-structure of $SU(m\ve n)$ induce on $\mf{g}^*$ a structure of graded real Lie superalgebra noted $(\mf{g}^*,\varphi_{G})$ (cf. prop $\ref{pp1}$) isomorph to $(\mf{b},\phi_{B})$. Similarly, the Poisson-Lie superbracket and the graded real $\bigstar$-structure of $s(AN)$ induce on $\mf{b}^*$ a structure of graded real Lie superalgebra noted $(\mf{b}^*,\varphi_{B})$ (cf. prop $\ref{pp1}$) isomorph to $(\mf{g},\phi_{G})$.
\label{th2}
\end{theorem}


\vskip1pc
\begin{rmk}
The theorem $(\ref{th2})$ states the duality of Poisson-Lie for the graded real supergroup $(SU(m\ve n),\bigstar)$ and $(s(AN),\bigstar)$, i.e. the duality with the graded $\bigstar$-structures.
\end{rmk}

\begin{rmk}
The superbracket $(\ref{crochetPL})$ have already been introduced in the article $\cite{Andru}$. Nevertheless, its used to define a double seems new.\\
\\
Moreover, this superbracket is basis independant. The operators used in the formula $(\ref{crochetPL})$ are defined by the following equalities
\ben
\n^{R}(f)=\d(f')f'', \quad \n^{L}(f)=(-1)^{\ve \n^{L} \ve (\ve f \ve +1)}f'\d(f''), \quad \forall f \in \slmn.
\label{isodef}
\een
These operators fullfil the super-Leibniz rule i.e.
\be
\n(fg)=\n(f)g+(-1)^{\ve \n \ve \ve f \ve}f\n(g),\quad  \forall f,g \in \slmn,
\ee
and the following properties
\be 
\Delta(\n^{R}(f))=(\n^{R}\o 1)\Delta (f), \; \Delta(\n^{L}(f))=(1\o \n^{L})\Delta (f), \; \forall f \in \slmn.
\ee
\end{rmk}



\vskip2pc
\no {\bf Proof of the Theorem $\ref{th1}$:}\\
\\
\underline{i) Poisson-Lie superbracket:}\\
\\
First, we show that $\Delta$ is a Poisson morphism (cf. prop. $iv)$ def. $\ref{defspl}$). With the help of the formula $(\ref{isodef})$ we rewrite the superbracket $(\ref{crochetPL})$ with the $\e$-derivations, i.e. we have
\be
\{f,g\}=(-1)^{\ve a \ve \ve g \ve }r^{ab}f'\, \d_{a}f''\, g'\, \d_{b}g''-(-1)^{\ve a \ve \ve f \ve}r^{ab}\d_{a}f' \,f''\, \d_{b}g' \,g''
\ee
where $r^{ab}=\frac{1}{2}(R_{\mf{d}}\hat{h}_{a},\hat{h}_{b})_{\mf{d}}$, $\d_{a}=\d_{h_{a}}$ and $f,g\in \slmn$.\\
On the one hand, we have
\begin{eqnarray}
\Delta\{f,g\}=(-1)^{\ve a \ve \ve g \ve + \ve g' \ve \ve f''\ve}r^{ab}f'g'\o f''g''\d_{a}f'''\d_{b}g'''\nonumber\\
-(-1)^{\ve a \ve \ve f \ve + \ve g'' \ve \ve f''' \ve}r^{ab}\d_{a}f'\d_{b}g'f''g''\o f'''g'''.
\label{pm1}
\end{eqnarray}
On the other hand, we have
\begin{eqnarray}
\{\Delta f,\Delta g\}&=&(-1)^{\ve f'' \ve \ve g' \ve}\{f',g'\}\o f''g''+(-1)^{\ve f'' \ve \ve g' \ve}f'g'\o \{f'',g''\}\nonumber\\
&=&(-1)^{\ve f''' \ve (\ve g' \ve+\ve g'' \ve)+a(\ve g' \ve + \ve g'' \ve)}r^{ab}f'g'\d_{a}f''\d_{b}g''\o f''' g'''\nonumber\\
&-&(-1)^{\ve f''' \ve (\ve g' \ve+\ve g'' \ve)+a(\ve f' \ve + \ve f'' \ve)}r^{ab}\d_{a}f'\d_{b}g'f''g''\o f'''g'''\nonumber\\
&+&(-1)^{\ve g' \ve (\ve f'' \ve+\ve f''' \ve)+a(\ve g'' \ve + \ve g''' \ve)}r^{ab}f'g'\o f''g''\d_{a}f'''\d_{b}g'''\nonumber\\
&-&(-1)^{\ve g' \ve (\ve f'' \ve+\ve f''' \ve)+a(\ve f'' \ve + \ve f''' \ve)}r^{ab}f'g'\o \d_{a}f''\d_{b}g''f'''g'''
\label{pm2}
\end{eqnarray}
We observe that $\d_{a}f$ vanish when $\ve f \ve \neq \ve a \ve$, which means that the a priori non-zero terms are such that $\ve f \ve=\ve a \ve$. Thus, the exposant of the sign of the first term becomes $(\ve f'''\ve +\ve f''\ve )(\ve g'\ve +\ve g''\ve )=(\ve f\ve +\ve f'\ve )(\ve g\ve +\ve g'''\ve )$ and the exposant of the sign of the fourth term is equal to $(\ve f''\ve +\ve f'''\ve )(\ve g''\ve +\ve g'\ve )=(\ve f\ve +\ve f'\ve )(\ve g\ve +\ve g'\ve )$. So that, the first and last term cancel each other. On the other hand, the exposant of the sign of the second term is $\ve a\ve \ve f\ve +\ve f''\ve \ve g''\ve $ and for the third one have $\ve f''\ve \ve g'\ve +\ve a\ve \ve g\ve $. Thus, the expression $(\ref{pm2})$ is equal to the expression $(\ref{pm1})$, which imply that $\Delta$ is a Poisson morphism for the superbracket $(\ref{crochetPL})$.\\ 
\\ We deduce the superantisymetry of the superbracket from the supersymetry of the scalar product $(.,.)_{\mf{d}}$ and the antisymetry of $R_{\mf{d}}$. This superbracket fullfils also the super-Leibniz rule because the operators $\n^{R,L}_{h_{a}}$ are superderivations.\\
\\ It remains to prove the super-Jacobi identity. We define the following Yang-Baxter superbracket\footnote[8]{This Yang-Baxter superbracket has been already define in $\cite{Andru}$.}
\be 
[\nabla^{L,R}_{r_{1}}\otimes \nabla^{L,R}_{r_{2}},\nabla^{L,R}_{r'}\otimes \nabla^{L,R}_{r''}]=[\nabla^{L,R}_{r_{1}},\nabla^{L,R}_{r'}] \otimes \nabla^{L,R}_{r''} \otimes \nabla^{L,R}_{r_{2}} + 
\ee
\be
\nabla^{L,R}_{r'} \otimes \nabla^{L,R}_{r_{1}} \otimes [\nabla^{L,R}_{r_{2}},\nabla^{L,R}_{r''}]-
\nabla^{L,R}_{r'} \otimes [\nabla^{L,R}_{r''},\nabla^{L,R}_{r_{1}}] \otimes \nabla^{L,R}_{r_{2}}.
\ee
We used also the following convention 
\be
 \nabla_{A} \otimes \nabla_{B} \otimes \nabla_{C}(f \otimes g \otimes h)= (-1)^{\ve f \ve (\ve B \ve + \ve C \ve)+ \ve g \ve \ve C \ve}
\nabla_{A}f \;  \nabla_{B}g \; \nabla_{C}h. 
\ee
Then, the super-Jacobi identity is equivalent to
\be
[\nabla^{L}_{r_{1}}\otimes \nabla^{L}_{r_{2}},\nabla^{L}_{r'}\otimes \nabla^{L}_{r''}]f \otimes g \otimes h +[\nabla^{R}_{r_{1}}\otimes \nabla^{R}_{r_{2}},\nabla^{R}_{r'}\otimes \nabla^{R}_{r''}]f \otimes g \otimes h=0,
\ee
where $r=\sum^{m+n}_{i,j=1}(R(\hat{h}_{i}),\hat{h}_{j}) h_{i} \otimes h_{j}=r' \otimes r''=r_{1}\otimes r_{2}$. Moreover, from the fact thar $R_{\mf{d}}$ fullfils the Yang-Baxter identity and is antisymetric we deduce that
\ben
[\nabla^{L}_{r_{1}}\otimes \nabla^{L}_{r_{2}},\nabla^{L}_{r'}\otimes \nabla^{L}_{r''}]=\sum^{m+n}_{i,j,k=1}(\hat{h}_{i},[\hat{h}_{j},\hat{h}_{k}])_{\mf{d}} \nabla^{L}_{h_{i}} \otimes \nabla^{L}_{h_{j}} \otimes \nabla^{L}_{h_{k}},
\label{somme1}
\een
\ben
[\nabla^{R}_{r_{1}}\otimes \nabla^{R}_{r_{2}},\nabla^{R}_{r'}\otimes \nabla^{R}_{r''}]=-\sum^{m+n}_{i,j,k=1}(\hat{h}_{i},[\hat{h}_{j},\hat{h}_{k}])_{\mf{d}} \nabla^{R}_{h_{i}} \otimes \nabla^{R}_{h_{j}} \otimes \nabla^{R}_{h_{k}}.
\label{somme2}
\een
These two expressions are basis independant, thereby we write them with the basis $e_{i}=(v_{i},0),(0,v_{i})$ and its dual basis $\hat{e}_{i}=(\hat{v}_{i},0),(0,-\hat{v}_{i})$. Furthermore, it is sufficient to prove the super-Jacobi identity on three generators taken among $u_{ij},w_{ab}$. Thus, we evaluate the terms of the sums $(\ref{somme1}),(\ref{somme2})$ on the generators $u_{ab},u_{cd},w_{mn}$, we have
\be
(-1)^{\ve u_{ab} \ve (\ve k \ve + \ve  j \ve)+\ve u_{cd} \ve \ve j \ve}((\hat{v}_{i},0),[(\hat{v}_{j} , 0),(0,-\hat{v}_{k})])_{\mf{d}} \nabla^{R,L}_{i}u_{ab} \;  \nabla^{R,L}_{j}u_{cd} \; \nabla^{R,L}_{k}v_{mn}.
\ee 
But the term $[(\hat{v}_{j} , 0),(0,-\hat{v}_{k})]$ is zero, this proves the super-Jacobi identity on three generators $u_{ab},u_{cd},w_{mn}$. The proof is the same for three generators $u_{ab},w_{cd},w_{mn}$. Therefore, it remains to prove the super-Jacobi identity on three generators of the same kind. We do it for three generators choosed among the generators $u_{ab}$, the proof is identical for the generators $w_{cd}$ thus we don't give it. We have the following equality (cf. annexe 2)
\begin{eqnarray}
\lefteqn{\sum^{m+n}_{i,j,k=1}(\hat{v}_{i},[\hat{v}_{j},\hat{v}_{k}])_{sl} \nabla^{R,L}_{v_{i}} \otimes \nabla^{R,L}_{v_{j}} \otimes \nabla^{R,L}_{v_{k}}={}} 
\nonumber\\ & & 
{} \sum_{a,b,c,d,s,t=1}^{m+n}(\hat{E}_{ab},[\hat{E}_{cd},\hat{E}_{st}])_{sl} \nabla^{R,L}_{E_{ab}} \otimes \nabla^{R,L}_{E_{cd}} \otimes \nabla^{R,L}_{E_{st}},  
\label{egalite}
\end{eqnarray}
where $\hat{E}_{st}=(-1)^{\ve  s \ve} 2i E_{ts}$ and fullfils  $(E_{pq},\hat{E}_{st})_{sl}= \delta_{ps}\delta_{qt}$. Then by a direct computation we find
\begin{eqnarray}
\lefteqn{ \sum_{a,b,c,d,s,t=1}^{m+n}(\hat{E}_{ab},[\hat{E}_{cd},\hat{E}_{st}])_{sl} \nabla^{L}_{E_{ab}}\otimes \nabla^{L}_{E_{cd}} \otimes \nabla^{L}_{E_{st}} (u_{ij}\otimes u_{kl} \otimes u_{pq})-{}} 
\nonumber\\ & & 
{} \sum_{a,b,c,d,s,t=1}^{m+n}(\hat{E}_{ab},[\hat{E}_{cd},\hat{E}_{st}])_{sl} \nabla^{R}_{E_{ab}} \otimes \nabla^{R}_{E_{cd}} \otimes \nabla^{R}_{E_{st}}(u_{ij}\otimes u_{kl} \otimes u_{pq})=0. \nonumber \\  
\label{egalite2}
\end{eqnarray}
Thus, we have proved the super-Jacobi identity.

\vskip2pc
\no \underline{ii) I,J Poisson ideals:}\\
\\
In view to prove that $I,J$ are Poisson ideals for the superbracket $(\ref{crochetPL})$, we exprime this superbracket in another basis. Let $T_{i},t_{i}$ be a basis of respectively $\mf{g}$ and $\mf{b}$ defined by (cf. annexe $1$ for notations)
\be
T_{i}=\{(E_{\alpha},E_{\alpha}),(E_{-\alpha},E_{-\alpha}),(H_{\mu},H_{\mu}),(\tilde{E}_{\beta},\tilde{E}_{\beta}),(\tilde{E}_{-\beta},\tilde{E}_{-\beta}),
\ee
\ben
(\tilde{H}_{\nu},\tilde{H}_{\nu}),(V_{\gamma},V_{\gamma}),(V_{-\gamma},V_{-\gamma}),(H_{0},H_{0})\} \;\;,
\label{base1}
\een
\be
t_{i}=\{(2iE_{-\alpha},0),(0,-2iE_{\alpha}),(iH_{\mu},-iH_{\mu}),(2i\tilde{E}_{-\beta},0),(0,-2i\tilde{E}_{\beta}),
\ee
\ben
(i\tilde{H}_{\nu},-i\tilde{H}_{\nu}),(2iV_{-\gamma},0),(0,2iV_{\gamma}),(iH_{0},-iH_{0})\}.
\label{base2}
\een
These basis fullfil $(T_{i},t_{j})_{\mf{d}}=\delta_{ij}$. Since $\mf{d}=\mf{g}\dotplus \mf{b}$, the vectors $T_{i},t_{i}$ are a basis of $\mf{d}$, and its dual basis is $\hat{T}_{i}=t_{i},\; \hat{t}_{i}=(-1)^{\ve i \ve}T_{i}$. As $R_{\mf{d}}=P_{\mf{b}}-P_{\mf{g}}$, we deduce that $R_{\mf{d}}(\hat{T}_{i})=t_{i}$ et $R_{\mf{d}}(\hat{t}_{i})=-(-1)^{\ve i \ve}T_{i}$. In this basis $T_{i},t_{i}$ the expression of the superbracket  $(\ref{crochetPL})$ becomes 
\begin{eqnarray}
\{ f, g \}& = &\frac{1}{2}[\sum_{a=1}^{n+m}(-1)^{\ve a \ve \ve f \ve}\nabla^{L}_{T_{a}}f \; \nabla^{L}_{t_{a}}g 
-(-1)^{\ve a \ve + \ve a \ve \ve f \ve }\nabla^{L}_{t_{a}}f \; \nabla^{L}_{T_{a}}g \nonumber \\ & &
-(-1)^{\ve a \ve \ve f \ve}\nabla^{R}_{T_{a}}f \; \nabla^{R}_{t_{a}}g 
+(-1)^{\ve a \ve + \ve a \ve \ve f \ve }\nabla^{R}_{t_{a}}f \; \nabla^{R}_{T_{a}}g],
\label{sb1}
\end{eqnarray}
Moreover, we have (cf. annexe 3)
\be
C^R=C^L
\ee
with
\be
C^{L,R}=\sum_{a=1}^{m+n}\nabla^{L,R}_{T_{a}}\otimes \nabla^{L,R}_{t_{a}}+(-1)^{\ve a \ve}\nabla^{L,R}_{t_{a}} \otimes \nabla^{L,R}_{T_{a}}.
\ee
From this equality, we deduce the following new expression of the superbracket  $(\ref{sb1})$ 
\ben
\{ f, g \}=\sum_{i=1}^{n+m}(-1)^{\ve i \ve \ve f \ve}(\nabla^{L}_{T_{i}}f \; \nabla^{L}_{t_{i}}g-\nabla^{R}_{T_{i}}f \; \nabla^{R}_{t_{i}}g),
\label{sb2}
\een
Then, since $T_{i}$ is a basis of $\mf{g}$, the ideal $I$ is invariant under the superderivations $\n^{R,L}_{T_{i}}$. Therefore, from the expression $(\ref{sb2})$ we have $\{ I,f\}\subset I$ for all $f\in \slmn$. Similarly, since $t_{i}$ is a basis of $\mf{b}$, $J$ is invariant under the superderivations  $\n^{R,L}_{t_{i}}$ and so we deduce from the expression $(\ref{sb2})$ that $\{f,J \}\subset J$ for all $f\in \slmn$. Thus, the Hopf ideal $I,J$ are Poisson ideals. Furthermore, it is clear that (because the superbracket $(\ref{crochetPL})$ is Poisson-Lie) the superbrackets on $\sumn$ and $\san$ are super Poisson-Lie.


\vskip2pc
\no \underline{iii) $\{f,g\}^*=\{f^*,g^*\}$:}\\
\\
First, we remark that
\be
\n^{R}_{M}(f^*)=(\n^{R}_{\phi(M)}f)^*,\quad \n^{L}_{M}(f^*)=(\n^{L}_{\phi(M)}f)^*.
\ee
Thus, we have
\begin{eqnarray*}
\{f^*,g^* \}&=&\sum_{a,b=1}^{m+n}(-1)^{\ve a \ve \ve f \ve}(R_{\mf{d}}(\hat{h}_{a}),\hat{h}_{b})_{\mf{d}}[\n^L_{h_{a}}f^* \; \n^L_{h_{b}}g^*-\n^R_{h_{a}}f^* \; \n^R_{h_{b}}g^*]\\
&=&\sum_{a,b=1}^{m+n}(-1)^{\ve a \ve \ve f \ve}(R_{\mf{d}}(\hat{h}_{a}),\hat{h}_{b})_{\mf{d}}[(\n^L_{\phi(h_{a})}f)^* \; (\n^L_{\phi(h_{b})}g)^*\\
&&-(\n^R_{\phi(h_{a})}f)^* \; (\n^R_{\phi(h_{b})}g)^*]
\end{eqnarray*}
In exchanging $h_{a}$ by $\phi(h_{a})$ and by observing that  $\widehat{\phi(h_{a})}=\phi(\hat{h}_{a})$, we deduce 
\begin{eqnarray*}
\{f^*,g^*\}&=&\sum_{a,b=1}^{m+n}(-1)^{\ve a \ve \ve f \ve}(R_{\mf{d}}(\widehat{\phi(h_{a})}),\widehat{\phi(h_{b})})_{\mf{d}}[(\n^L_{\phi^2(h_{a})}f)^* \; (\n^L_{\phi^2(h_{b})}g)^*\\
&& -(\n^R_{\phi^2(h_{a})}f)^* \; (\n^R_{\phi^2(h_{b})}g)^*]\\
&=&\sum_{a,b=1}^{m+n}(-1)^{\ve a \ve \ve f \ve}\ov{(R_{\mf{d}}(\hat{h}_{a}),\hat{h}_{b})_{\mf{d}}}[(\n^L_{h_{a}}f)^* \; (\n^L_{h_{b}}g)^*-(\n^R_{h_{a}}f)^* \; (\n^R_{h_{b}}g)^*]\\
&=&\{f,g \}^*.
\end{eqnarray*}

\vskip2pc
\no{\bf Proof of the corollary $\ref{co1}$:}\\
\\
It is a direct consequence of the previous theorem $\ref{th1}$.\hfill$\blacksquare$


\vskip2pc
\no{\bf Proof of the theorem $\ref{th2}$:}\\
\\
The strategy of the proof is the following, first we prove that the Poisson-Lie superbracket $(\ref{sb2})$ and the $*$-structure of $\slmn$ induce (because of the prop. $\ref{pp1}$) respectively a structure of direct sum of Lie superalgebra on  $\mf{d}^*=\m{G}^* \oplus \m{B}^*$ and a graded real structure $\varphi$ on $\mf{d}^*$, which leave invariant  $\m{G}^*$ and $\m{B}^*$. Here, $\m{G}^*$  ($\m{B}^*$) is defined  like the space of linear forms on $\mf{d}=\mf{g} \dotplus \mf{b}$ which cancel on $\mf{b}$ ($\mf{g}$) i.e. $\m{G}^*$ and $\m{B}^*$ are naturally identified to $\mf{g}^*$ and $\mf{b}^*$ (i.e. with the linear dual of $\mf{g}$ and $\mf{b}$). Next, we show that $(\m{G}^*,\pg)$ and $(\m{B}^*,\pb)$ are respectively isomorph to $(\mf{b},\phi_{B})$ and $(\mf{g},\phi_{G})$. On the other hand, the dual space of $\mf{g}$ and $\mf{b}$ noted $\mf{g}^*,\, \mf{b}^*$ are also defined directly from  $\sumn$ and $\san$ as \footnote[9]{We frequently note in this proof the maps on $\sumn$ ($\san$) which come form maps on $\slmn$ by the index $G$ ($B$). Thus, $\e_{G}$ is the counity of $\sumn$. Similarly, we index by the letter $D$ the maps defined on $\slmn$, for instance $\e_{D}$ is the counity of $\slmn$.} $\mf{g}^*=Ker\e_{G}/(Ker\e_{G})^2$ and $\mf{b}^*=Ker\e_{B}/(Ker\e_{B})^2$. Then, since $\sumn$
and $\san$ are graded real Poisson-Lie supergroups, $\mf{g}^*$ and $\mf{b}^*$ are graded real Lie superalgebra which we note respectively $(\mf{g}^*,\varphi_{G})$ and $(\mf{b}^*,\varphi_{B})$. Finaly, we show that $(\mf{g}^*,\varphi_{G})$ and $(\mf{b}^*,\varphi_{B})$  are isomorph respectively to $(\m{G}^*,\pg)$ and $(\m{B}^*,\pb)$. To conclude, we have proved that $(\mf{g}^*,\varphi_{G})$ and $(\mf{b}^*,\varphi_{B})$ are isomorph respectively to ($\mf{b}$,$\phi_{B}$) and ($\mf{g}$,$\phi_{G}$). These show that $\sumn$ and $\san$ are duals Poisson-Lie supergroup. \\
\\

\no Let $f_{i},g_{i}$ be two sets of elements of $\slmn$ with $i=1...(m+n)^2-1$ defined by
\be
f_{i}=\{\sum_{i,j=1}^{m+n}(E_{\alpha})_{ij}u_{ij},\,\sum_{i,j=1}^{m+n}(E_{-\alpha})_{ij}w_{ij},\,\sum_{i=1}^{m+n}-\frac{(H_{\mu})_{ii}}{2}(S(u_{ii})-w_{ii}),
\ee
\be
\sum_{i,j=1}^{m+n}(\tilde{E}_{\beta})_{ij}u_{ij},\,\sum_{i,j=1}^{m+n}(\tilde{E}_{-\beta})_{ij}w_{ij},\,\sum_{i=1}^{m+n}-\frac{(\tilde{H}_{\nu})_{ii}}{2}(S(u_{ii})-w_{ii}),
\ee
\ben
\sum_{i,j=1}^{m+n}(V_{\gamma})_{ij}u_{ij},\,\sum_{i,j=1}^{m+n}(V_{-\gamma})_{ij}w_{ij},\,\sum_{i=1}^{m+n}-\frac{(H^*_{0})_{ii}}{2}(S(u_{ii})-w_{ii})\},
\label{fctf}
\een
\ben
g_{i}=\{\sum_{k,l}^{m+n}((A_{i})_{kl}u_{kl}+(B_{i})_{kl}w_{kl})\},
\label{fctg}
\een
with
\be
\{(A_{i},B_{i})\}=\{(-\frac{i}{2}E_{\alpha},\frac{i}{2}E_{\alpha}), (-\frac{i}{2}E_{-\alpha},\frac{i}{2}E_{-\alpha}),(-\frac{i}{2}H_{\mu},\frac{i}{2}H_{\mu}),(\frac{i}{2}\tilde{E}_{\beta},-\frac{i}{2}\tilde{E}_{\beta}),
\ee
\be
(\frac{i}{2}\tilde{E}_{-\beta},-\frac{i}{2}\tilde{E}_{-\beta}),(\frac{i}{2}\tilde{H}_{\nu},-\frac{i}{2}\tilde{H}_{\nu}),(-\frac{i}{2}V_{\gamma},\frac{i}{2}V_{\gamma}),(\frac{i}{2}V_{-\gamma},-\frac{i}{2}V_{-\gamma}),(-\frac{i}{2}H^*_{0},\frac{i}{2}H^*_{0}) \},
\ee
where $H_{0}^*=\frac{n-m}{n+m}H_{0}$. We remark that $f_{i}\in J$ and $g_{i}\in I$. Moreover, we note $\O_{D}$ the canonical projection
\be
\O_{D}:Ker\e_{D}\to Ker\e_{D}/(Ker\e_{D})^2.
\ee
For $\d_{(A,B)}\in \mf{d}$ and $f \in Ker\e_{D}$, the duality between $\mf{d}$ and $\mf{d}^*=Ker\e_{D}/(Ker\e_{D})^2$ is given by (cf. $(\ref{dualite})$) 
\be
\langle \d_{(A,B)},\O_{D}(f) \rangle=\d_{(A,B)}(f).
\ee
Thus the familly of vectors $\O_{D}(f_{i})$ are in duality with $\d_{T_{i}}$ and cancel on $\d_{t_{i}}$, while $\O_{D}(g_{i})$ are in duality with $\d_{t_{i}}$ and vanish on $\d_{T_{i}}$. So, $\O_{D}(f_{i})$ and $\O_{D}(g_{i})$ are respectively a basis of $\m{G}^*$ and $\m{B}^*$. From the proposition $\ref{pp1}$, we deduce that $\mf{d}^*$ is endowed with the following Lie superbracket
\ben
[\O_{D}(f),\O_{D}(g)]_{\mf{d}^*}=\O_{D}(\{ f,g\})
\label{crochetlie}
\een
for all $f,g \in Ker\e_{D}$. Now,  we show that $\mf{d}^*$ is a direct sum of Lie superalgebra i.e. $\mf{d}^*=\m{G}^*\oplus\m{B}^*$. 
The expression of the Lie superbracket $(\ref{crochetlie})$ in the basis $\O_{D}(f_{i}),\O_{D}(g_{i})$ is 
\be
[\O_{D}(f),\O_{D}(g)]_{\mf{d}^*}=\sum_{k=1}^{(m+n)^2-1}\d_{T_{k}}(\{ f,g\})\O_{D}(f_{k})+\sum_{s=1}^{(m+n)^2-1}\d_{t_{s}}(\{ f,g\})\O_{D}(g_{s}),
\ee
for all $f,g\in \slmn$. Since $\d_{T_{k}}$ and $\d_{t_{s}}$ cancel respectively on $I,J$ and that $\{f_{i},g_{j}\}$  is in $I$ and $J$, we deduce that 
\be
[\O_{D}(f_{i}),\O_{D}(g_{j})]_{\mf{d}^*}=0.
\ee
We have also
\ben
[\O_{D}(f_{i}),\O_{D}(f_{j})]_{\mf{d}^*}=\sum_{k=1}^{(m+n)^2-1}\d_{T_{k}}(\{ f_{i},f_{j}\})\O_{D}(f_{k}),
\label{salgg}
\een
\ben
[\O_{D}(g_{i}),\O_{D}(g_{j})]_{\mf{d}^*}=\sum_{k=1}^{(m+n)^2-1}\d_{t_{k}}(\{ g_{i},g_{j}\})\O_{D}(g_{k}).
\label{salgb}
\een
We find with the help of the expression $(\ref{sb2})$ of the superbracket $(\ref{crochetPL})$ the equalities 
\be
\d_{T_{k}}(\{ f_{i},f_{j}\})=(-1)^{\ve i \ve +1}[\d_{t_{i}},\d_{T_{k}}](f_{j})
,\; \d_{t_{k}}(\{ g_{i},g_{j}\})=(-1)^{\ve i \ve \ve j \ve+1}[\d_{T_{j}},\d_{t_{k}}](g_{i}).
\ee  
Moreover, from the following equality
\be
[\d_{T_{a}},\d_{t_{b}}]=\sum_{k=1}^{(m+n)^2-1}([T_{a},t_{b}],t_{k})_{\mf{d}}\,\d_{T_{k}}+\sum_{s=1}^{(m+n)^2-1}(T_{k},[T_{a},t_{b}])_{\mf{d}}\,\d_{t_{k}}
\ee
we deduce
\ben
[\O_{D}(f_{i}),\O_{D}(f_{j})]_{\mf{d}^*}=\sum_{k=1}^{(m+n)^2-1}(-1)^{\ve i \ve \ve  j \ve}(T_{k},[t_{i},t_{j}])_{\mf{d}}\,\O_{D}(f_{k}),
\label{eqcr1}
\een
\ben
[\O_{D}(g_{i}),\O_{D}(g_{j})]_{\mf{d}^*}=\sum_{k=1}^{(m+n)^2-1}(-1)^{\ve i \ve \ve j \ve+1}([T_{i},T_{j}],t_{k})_{\mf{d}}\,\O_{D}(g_{k}).
\label{eqcr2}
\een
Thus, $\mf{d}^*$ is a direct sum of Lie superalgebra.
On the other hand, with the help of the proposition $\ref{pp1}$, it turns out that $\mf{d}^*$ is endowed with the graded real struture
\be
\varphi(\O_{D}(f))=\O_{D}(f^\bigstar),
\ee
for all $f \in \slmn$. Furthermore, as $I$ and $J$ are invariant under $\bigstar$ and $f_{i}$, $g_{j}$ are respectively in $I$ and $J$, it turns out that $\m{G}^*$ and $\m{B}^*$ are invariant under $\varphi$. So $(\m{G}^*,\pg)$ and $(\m{B}^*,\pb)$ are graded real Lie superalgebra.\\
\\
Now, we show that $(\m{G}^*,\pg)$ is isomorph to $(\mf{b},\phi_{B})$. Let $\O_{D}(\mf{f}_{a})=i^{\ve a \ve}\O_{D}(f_{a})$. We have the following equality
\be
[\O_{D}(\mf{f}_{a}),\O_{D}(\mf{f}_{b})]_{\mf{d}^*}=\sum_{c=1}^{(m+n)^2-1}(T_{k},[t_{i},t_{j}])_{\mf{d}}\,\O_{D}(\mf{f}_{k}).
\label{eqqq1}
\ee
Let $\m{T}$ be the isomorphism defined by 
\be
\m{T}:\O_{D}(\mf{f}_{a})\to \d_{t_{a}}
\ee
then the equality $(\ref{eqcr1})$ imply
\be
\m{T}([\O_{D}(\mf{f}_{a}),\O_{D}(\mf{f}_{b})]_{\mf{d}^*})=[\m{T}(\O_{D}(\mf{f}_{a})),\m{T}(\O_{D}(\mf{f}_{b}))]_{\mf{d}}
\ee
in other words  $\m{T}$ is an isomorphism of Lie superalgebra.
Nex with the help of the scalar product $(.,.)_{\mf{d}}$ we define the isomorphism $\mf{I}:\mf{d}^*\to \mf{d}$ by 
\be
\langle x, \omega \rangle=(x,\mf{I}(w))_{\mf{d}},
\ee 
where $x\in \mf{d}$ and $\omega \in \mf{d}^*$. Thus, $\mf{I}(\O_{D}(\mf{f}_{a}))=i^{\ve a \ve}\d_{t_{a}}$. We remark that $\mf{I}(\omega)=i^{\ve \omega \ve}\m{T}(\omega),\;\forall \omega \in \m{G}^*$. Moreover, we have
\be
\varphi(\O_{D}(f))=\sum_{a=1}^{(m+n)^2-1}\d_{T_{a}}(f^*)\,\O_{D}(f_{a}).
\ee 
We apply $\m{T}$ on the two members of the previous equality
\begin{eqnarray*}
\m{T}(\varphi(\O_{D}(f)))&=&\sum_{a=1}^{(m+n)^2-1}\d_{T_{a}}(f^*)\,\m{T}(\O_{D}(f_{a}))\\
&=&\sum_{a=1}^{(m+n)^2-1}(-i)^{\ve a \ve}\overline{\langle \phi(\d_{T_{a}}) , \O_{D}(f) \rangle}\,\d_{t_{a}}\\
&=&\sum_{a=1}^{(m+n)^2-1}(-i)^{\ve a \ve}\overline{(\phi(\d_{T_{a}}),\mf{I}(\O_{D}(f)))}_{\mf{d}}\,\d_{t_{a}}\\
&=&\sum_{a=1}^{(m+n)^2-1}(i)^{\ve a \ve} (\d_{T_{a}},\phi(\mf{I}(\O_{D}(f))))_{\mf{d}}\,\d_{t_{a}}\\
&=&i^{\ve f \ve}\,\phi(\mf{I}(\O_{D}(f))).
\end{eqnarray*}
From the fact that $\mf{I}(\omega)=i^{\ve \omega \ve}\m{T}(\omega),\;\forall \omega \in \m{G}^*$ and the previous equality we find
\be
\m{T}(\varphi(\O_{D}(f)))=\phi(\m{T}(\O_{D}(f)))
\ee
i.e.
\be
\m{T}\circ \varphi \circ \m{T}^{-1}=\phi.
\ee
Thus, $(\m{G}^*,\pg)$ et $(\mf{b},\phi_{B})$ are graded real Lie superalgebra.\\
\\
We turn to the proof of the isomorphism of $(\m{B}^*,\pb)$ and $(\mf{g},\phi_{G})$. Let $\O_{D}(\mf{g}_{a})=(-1)^{\ve a \ve + 1}\O_{D}(g_{a})$ and the morphism $\m{S}:\m{B}^*\to\mf{g}$ define by 
\be
\m{S}(\O_{D}(\mf{g}_{a}))=\d_{T_{a}}.
\ee
Thus, The equality $(\ref{eqcr2})$ imply 
\be
\m{S}([\O_{D}(\mf{g}_{a}),\O_{D}(\mf{g}_{b})]_{\mf{d}^*})=[\m{S}(\O_{D}(\mf{g}_{a})),\m{S}(\O_{D}(\mf{g}_{b}))]_{\mf{d}}
\ee
i.e. $\m{S}$ is an isomorphism of Lie superalgebra from $\m{B}^*$ to $\mf{g}$. The scalar product gives also an isomorphism from $\mf{B}^*$ to $\mf{g}$ such that $\mf{I}(\O_{D}(\mf{g}_{a}))=-\d_{T_{a}}$. It appears that $\mf{I}=-\m{S}$. On the other hand, the graded real structure is defined on $\mf{B}^*$ by 
\be
\varphi(\O_{D}(g))=\O_{D}(g^*).
\ee 
We exprim this previous element in the basis $\O_{D}(g_{a})$ i.e.
\be
\varphi(\O_{D}(g))=\sum_{a=1}^{(m+n)^2-1}\overline{\langle \phi(\d_{t_{a}}), \O_{D}(g)\rangle} \O_{D}(g_{a}).
\ee 
We apply $\m{S}$ on each side of the previous equality, thus we have
\begin{eqnarray*}
\m{S}(\varphi(\O_{D}(g)))&=&\sum_{a=1}^{(m+n)^2-1}\overline{(\phi(\d_{t_{a}},\mf{I}(\O_{D}(g))))}_{\mf{d}}\m{S}(\O_{D}(g_{a}))\\
&=&\sum_{a=1}^{(m+n)^2-1}(-1)^{\ve a \ve+1}(\phi^2(\d_{t_{a}},\phi(\mf{I}(\O_{D}(g)))))_{\mf{d}}\d_{T_{a}}\\
&=&-(-1)^{\ve g \ve}\sum_{a=1}^{(m+n)^2-1}(\phi(\mf{I}(\O_{D}(g)))),\d_{t_{a}})_{\mf{d}}\d_{T_{a}}\\
&=&-(-1)^{\ve g \ve }\phi(\mf{I}(\O_{D}(g))).
\end{eqnarray*}
Thus, if $\ve g \ve=0$ we have $\m{S}\circ \varphi=-\phi \circ \mf{I}$, but since $\mf{I}=-\m{S}$ we deduce $\m{S}\circ \varphi \circ \m{S}^{-1}=\phi$. While for $\ve g \ve=1$ we have $\m{S}\circ \varphi=\phi\circ \mf{I}$ i.e. $\m{S}\circ \varphi \circ \m{S}^{-1}=-\phi$. In fact, we have 
\be
\m{S}\circ \varphi \circ \m{S}^{-1}(M)=(-1)^{\ve M \ve}\phi(M)=K\phi(K^{-1}M)
\ee 
with $M=(A,B),\; A,B \in sl(m\ve n,\C)$, $K=diag(1_{m},-1_{n})$ and $KM=(KA,KB)$. Thus, $\m{S}\circ \varphi \circ \m{S}^{-1}$ is conjugated to $\phi$ and so $(\m{B}^*,\pb)$ and $(\mf{g},\phi_{G})$ are isomorphic graded real Lie superalgebras.
\vskip3pc
\no Now,  we equipp $\mf{g}^*$ with a graded real Lie superalgebra, and we prove that it is isomorph to $(\m{G}^*,\pg)$. 
Let $i_{G},\O_{G}$ be the following morphism of superalgebra
\be
i_{G}:\slmn\to \sumn, \quad \O_{G}:Ker\e_{G}\to Ker\e_{G}/(Ker\e_{G})^2, 
\ee
and remember that $\mf{g}^*=Ker\e_{G}/(Ker\e_{G})^2$ by definition. The Lie superalgebra of $\sumn$ is $\mf{g}=\{\e{\rm-derivation}\; \d:\sumn \to \C \}$, but every  $\e$-derivations of $\mf{g}$ come from an $\e$-derivation of $\slmn$ which cancel on $I$. Thus, $(\d_{T_{i}})$ is an element of $\mf{g}$ defined by the following formula
\be
(\d_{T_{i}})(i_{G}(f))=\d_{T_{i}}(f)
\ee
with $f \in \slmn$ and $\d_{T_{i}}\in \mf{d}$ such that $\d_{T_{i}}(I)=0$.
Let  $\O_{G}(i_{G}(f_{i}))$ be vectors of $\mf{d}^*$, where $f_{i}$ are defined by  $(\ref{fctf})$, it is in duality with the basis $(\d_{T_{i}})$ of $\mf{d}$, indeed we have
\be
\langle (\d_{T_{i}}), \O_{G}(i_{G}(f_{j})) \rangle =(\d_{T_{i}})(i_{G}(f_{j}))=\d_{T_{i}}(f_{j})=\d_{ij}.
\ee 
So, the set $\O_{G}(i_{G}(f_{i}))$ is a basis of $\mf{g}^*$.
We define the Lie superbarcket by
\ben
[\O_{G}(i_{G}(f)),\O_{G}(i_{G}(g))]_{\mf{g}^*}=\O_{G}(i_{G}(\{f,g\}))
\label{defsalglie}
\een
for all $f,g \in Ker\e_{D}$.
On the other hand, we deduce from the equation $(\ref{salgg})$ that 
\be
\{ f_{i},f_{j}\}=\sum_{k=1}^{(m+n)^2-1}\d_{T_{k}}(\{f_{i},f_{j}\})f_{k}+(Ker\e_{D})^2.
\ee
We evaluate $\O_{G}\circ i_{G}$ on each side of the previous equation i.e.
\be
\O_{G}(i_{G}(\{f_{i},f_{j}\}))=\sum_{k=1}^{(m+n)^2-1}\d_{T_{k}}(\{f_{i},f_{j}\})\O_{G}(i_{G}(f_{k})).
\ee
Thus, the isomorphism of supervector space $\mf{T}:\mf{g}^*\to \m{G}^*$ defined by $\mf{T}(\O_{G}(i_{G}(f_{i})))=\O_{D}(f_{i})$ is an isomorphism of Lie superalgebra since from the previous equality and the formula $(\ref{salgg})$, $(\ref{defsalglie})$ we deduce 
\be
\mf{T}([\O_{G}(i_{G}(f)),\O_{G}(i_{G}(g))]_{\mf{g}^*})=[\mf{T}(\O_{G}(i_{G}(f))),\mf{T}(\O_{G}(i_{G}(g)))]_{\m{G}^*}.
\ee
Finaly, the graded real structure of $\mf{g}^*$ is defined by 
\be
\varphi_{G}(\O_{G}(i_{G}(f)))=\O_{G}(i_{G}(f^*)).
\ee
But, from the expression of $\O_{G}(i_{G}(f^*))$ in the basis $\O_{G}(i_{G}(f_{i}))$, we find
\begin{eqnarray*}
\varphi_{G}(\O_{G}(i_{G}(f^*)))&=&\sum_{i=1}^{(m+n)^2-1}\langle (\d_{T_{i}}),\O_{G}(i_{G}(f^*))\rangle \O_{G}(i_{G}(f_{i}))\\
&=&\sum_{i=1}^{(m+n)^2-1}\d_{T_{i}}(f^*)\O_{G}(i_{G}(f_{i})).
\end{eqnarray*}
We apply the morphism $\mf{T}$ on the previous equality
\begin{eqnarray*}
\mf{T}(\varphi_{G}(\O_{G}(i_{G}(f))))&=&\sum_{i=1}^{(m+n)^2-1}\d_{T_{i}}(f^*)\O_{D}(f_{i})\\
&=&\langle \d_{T_{i}},\O_{D}(f^*)\rangle \O_{D}(f_{i})\\
&=&\O_{D}(f^*)\\
&=&\pg(\O_{D}(f)).
\end{eqnarray*}
This imply 
\be
\mf{T}\circ \varphi_{G}\circ \mf{T}^{-1}=\pg. 
\ee
In other words $\mf{T}$ is a graded real Lie superalgebra morphism.\\
\\
In the same way, we can prove that $\mf{b}^*$ is isomorph like graded real Lie superalgebra to  $(\m{B}^*,\varphi)$, the proof is similar thus we don't write it. 
\hfill$\blacksquare$
\vskip3pc
\no Now, we explain why this construction is the super-version of the Lu-Weinstein Drinfeld double.\\
\\
For $n=0$, we have $\mf{d}^{\phi}=\mf{g}^{\phi}\dotplus \mf{b}^{\phi}=sl(m,\C)^{\R}=su(m)\dotplus an$ where $\mf{d}^{\phi}=\{ x \in \mf{d}: \phi(x)=x\}$. Since the fixed points of  $\phi$ is the set $\{(X,-\bar{X}^t),\; X \in sl(m,\C)\}$, the scalar product $(.,.)_{\mf{d}}$ becomes
\be
((X,-\bar{X}^t),(Y,-\bar{Y}^t))_{\mf{d}}=Im(tr(XY)).
\ee
Thus, what we have done is the complexified version of the Lu-Weinstein Drinfeld double, the fact that this complexified double gives (for $n=0$) a real double comes from the properties of  the real structure $\phi$ with the structure of Lie algebra, the scalar product $(.,.)_{\mf{d}}$ and the classical $r$-matrix  $R_{\mf{d}}$.


\subsection{Annexe 1: $sl(m\ve n,\C)$.}

In this annexe, we give some useful properties for the Lie superalgebra $sl(m\ve n,\C)$ for $m\neq n$. We assume in this section  $m\neq n$.\\
\\
The Lie superalgebra $\sl=\{M\in M(m\ve n,\C )\; / Str(M)=0\}$ have the triangular decomposition
\be
\sl=\sl_{+}\oplus \sl_{0} \oplus \sl_{-},
\ee 
where $\sl_{+}$ is the strict (zero on the diagonal) upper triangular complex matrices, $\sl_{0}$ is the set of supertraceless diagonal supermatrices and $\sl_{-}$ is the set of strict lower triangular supermatrices. Thus, every elements $x\in \sl$ have the unique decomposition $x=x_{+}+x_{0}+x_{-}$ where $x_{i}\in \sl_{i},\; i \in \{+,0,-\}$. On the other hand,  $sl(m\ve n,\C)$ is provided with the following scalar product
\ben
(M,N)_{sl}=-\frac{i}{2}Str(MN),\quad \forall M,N \in \sl,
\label{prodsc}
\een  
which satisfies the following properties
\be
([M,N],P)_{sl}=(M,[N,P])_{sl},\quad (M,N)_{sl}=(-1)^{\ve M \ve \ve N \ve}(N,M)_{sl},  
\ee
$\forall M,N,P \in \sl$. When a scalar product fullfils the first property we said that it is invariant, for the second property we said that it is supersymetric. Furthermore, $\sl$ has the following properties:\\
\\
$(1)$ $\sl_{i}$ are Lie sub-superalgebras for $i\in \{+,0,-\}$ and $[\sl_{0},\sl_{0}\oplus \sl_{\pm}]\subset \sl_{\pm}$,\\
\\
$(2)$ $\sl_{+}\subset \sl_{+}^{\bot},\quad \sl_{-}\subset \sl_{-}^{\bot},\quad \sl_{0}\in (\sl_{0}\oplus \sl_{-})^{\bot}$ where $g^{\bot}$ is the orthogonal of the set $g$ in  $\sl$ for the scalar product $(\ref{prodsc})$.\\
\\
Now, we give a basis of $sl(m\ve n,\C)$. We note $E_{st}$ the supermatrice such that $(E_{st})_{ij}=\d_{si}\d_{tj}$. A basis of $\sl_{+}$ is
\be
E_{\alpha}=E_{st} \;\;\;  1\leq s < t \leq m \; ,\; \tilde{E}_{\beta}=iE_{st} \;\;\; m+1\leq s<t\leq m+n  ,   
\ee
\be
V_{\gamma}=E_{st} \;\;\; m+1\leq t \leq m+n \;\; , \;\; 1\leq s\leq m.
\ee
these are the positive roots of $\sl$, which we symbolize by $\alpha,\gamma>0$. A basis of $\sl_{-}$ is
\be
E_{-\alpha}=E_{\alpha}^{t} \;\; , \;\; \tilde{E}_{-\beta}=\tilde{E}_{\beta}^{t} \;\; , \;\; V_{-\gamma}=V_{\gamma}^{t}  ,
\ee
where $t$ is the usual transposition. These are the negatives roots $\alpha,\beta,\gamma<0$. While a basis of $\sl_{0}$ is 
\be
H_{k}=\frac{1}{\sqrt{k(k+1)}}(\sum_{l=1}^{k}E_{ll}-kE_{k+1 k+1}), \;\;  k=1...m-1  , 
\ee 
\be
\tilde{H}_{s}=\frac{i}{\sqrt{s(s+1)}}(\sum_{l=m+1}^{s}E_{ll}-sE_{s+1 s+1})\; ,\;\;  s=m+1...m+n-1 ,
\ee
to which we must add the following supervector   
\be
H_{0}=\frac{\sqrt{n}}{\sqrt{m(n-m)}}(\sum_{k=1}^{m}E_{kk}+\frac{m}{n}        \sum_{k=m+1}^{m+n}E_{kk}),\; n>m
\ee
or 
\be
H_{0}=\frac{i\sqrt{n}}{\sqrt{m(m-n)}}(\sum_{k=1}^{m}E_{kk}+\frac{m}{n}        \sum_{k=m+1}^{m+n}E_{kk}),\; n<m
\ee
Thus a basis of  $\sl$ is  $v_{i}=(E_{\alpha},H_{\mu},\tilde{E}_{\beta},\tilde{H}_{\nu},V_{\gamma},H_{0}),\; i=1...m+n-1$. This basis is normalised such that
\be
Str(E_{\alpha}E_{-\beta})=\delta_{\alpha \beta},\;\;\; 
Str(\tilde{E}_{\alpha}\tilde{E}_{-\beta})=\delta_{\alpha \beta},\;\;\;
Str(V_{\gamma}V_{-\epsilon})=\delta_{\gamma \epsilon},
\ee
\be
Str(H_{\mu}H_{\nu})=\delta_{\mu \nu},\;\;\; Str(\tilde{H}_{\mu}\tilde{H}_{\nu})=\delta_{\mu \nu},\;\;\; Str(H_{0}H_{0})=1,
\ee
and the supertrace of other couples is zero. Furthermore, a basis $\hat{v}_{i}$ is a dual basis of $v_{i}$ if it fullfils $(v_{i},\hat{v}_{j})_{sl}=\d_{ij}$. Thus, the dual basis of  $v_{i}$ is $\hat{v}_{i}=(\hat{E}_{\alpha},\hat{H}_{\mu},\hat{\tilde{E}}_{\beta},\hat{\tilde{H}}_{\nu},\hat{V}_{\gamma},\hat{H}_{0})$ such that
\be
\hat{E}_{\alpha}=2iE_{-\alpha},\quad \hat{H}_{\mu}=2iH_{\mu},\quad 
\hat{\tilde{E}}_{\beta}=2i\tilde{E}_{-\beta}, \quad
\hat{\tilde{H}}_{\nu}=2i\tilde{H}_{\nu},
\ee
\be
\hat{V}_{\gamma}=2iV_{-\gamma},\quad \hat{V}_{-\gamma}=-2iV_{\gamma}, \quad
\hat{H}_{0}=2iH_{0},
\ee
with $\gamma>0$ and $\alpha,\beta$  positive as well as negative roots. 


\subsection{Annexe 2: Proof of the equality $(\ref{egalite})$.}

\no Let $e_{i}=\{E_{\alpha},H_{\mu}\}$ be a basis of $sl(m,\C)\subset sl(m\ve n,\C)$ and $\hat{e}_{i}=\{\hat{E}_{\alpha},\hat{H}_{\mu} \}$ the dual basis for $(.,.)_{sl}$, then we define the Wess-Zumino-Witten form 
\begin{eqnarray*}
WZW_{sl(m,\C)}&=&\sum_{i,j,k=1}^{m}(\hat{e}_{i},[\hat{e}_{j},\hat{e}_{k}])_{sl}\n^{R,L}_{e_{i}}\o \n^{R,L}_{e_{j}}\o\n^{R,L}_{e_{k}}\\
&=&\frac{1}{6}\sum_{i,j,k=1}^{m}(\hat{e}_{i},[\hat{e}_{j},\hat{e}_{k}])_{sl}\n^{R,L}_{e_{i}}\w \n^{R,L}_{e_{j}}\w\n^{R,L}_{e_{k}}. 
\end{eqnarray*}
We develop the right hand side of this equation and gather identical terms i.e. we obtain 
\begin{eqnarray*}
WZW_{sl(m,\C)} & = & \frac{(2i)^{3}}{6}(\sum_{\a,\b,\g}(E_{-\a},[E_{-\b},E_{-\g}])_{sl}\n^{R,L}_{E_{\a}}\w \n^{R,L}_{E_{\b}}\w \n^{R,L}_{E_{\g}}\\
& & +6\sum_{\a>0,\b>0,\mu}\n^{R,L}_{E_{-\a}}\w\n^{R,L}_{E_{\b}}\w\n^{R,L}_{[E_{\a},E_{-\b}]})
\end{eqnarray*}
On the other hand, we have the following equality
\be
\sum_{s,t,p,q,a,b=1}^{m}(E_{st},[E_{pq},E_{ab}])_{sl} \n^{R,L}_{E_{st}}\w \n^{R,L}_{E_{pq}} \w \n^{R,L}_{E_{ab}}=
\ee
\be
\sum_{s\neq t,p\neq q, a\neq b=1}^{m}(E_{ts},[E_{qp},E_{ba}])_{sl}\, \n^{R,L}_{E_{st}}\w \n^{R,L}_{E_{pq}} \w \n^{R,L}_{E_{ab}}+\, 6 \sum_{s>t,p>q=1}^{m} \n^{R,L}_{E_{st}}\w \n^{R,L}_{E_{pq}} \w \n^{R,L}_{[E_{ts},E_{qp}]}
\ee
Since $E_{\a}=E_{st}$ with $t>s$ for $\a >0$ and $s>t$ for $\a<0$, we deduce  
\be
WZW_{sl(m,\C)}=\frac{(2i)^{3}}{6}\sum_{s,t,p,q,a,b=1}^{m}(E_{st},[E_{pq},E_{ab}])_{sl} \n^{R,L}_{E_{st}}\w \n^{R,L}_{E_{pq}} \w \n^{R,L}_{E_{ab}}.
\ee
or  
\ben
WZW_{sl(m,\C)}=\frac{1}{6}\sum_{s,t,p,q,a,b=1}^{m}(\hat{E}_{ts},[\hat{E}_{qp},\hat{E}_{ba}])_{sl} \n^{R,L}_{E_{st}}\w \n^{R,L}_{E_{pq}} \w \n^{R,L}_{E_{ab}}.
\label{wzw1}
\een
with $\hat{E}_{st}=(-1)^{\ve s\ve}2iE_{ts}$ such that $(E_{ts},\hat{E}_{ab})=\d_{t,b}\d_{s,a}$.
In the same way, we prove that 
\ben
WZW_{sl(n,\C)}=
\frac{1}{6}\sum_{s,t,p,q,a,b=m+1}^{m+n}(\hat{E}_{ts},[\hat{E}_{qp},\hat{E}_{ba}])_{sl} \n^{R,L}_{E_{st}}\w \n^{R,L}_{E_{pq}} \w \n^{R,L}_{E_{ab}}.
\label{wzw2}
\een
where we have define $WZW_{sl(n,\C)}$ in the basis $\tilde{e}_{i}=\{\tilde{E}_{\a},\tilde{H}_{\nu}\}$ and dual basis $\hat{\tilde{e}}_{i}=\{\hat{\tilde{E}}_{\a},\hat{\tilde{H}}_{\nu}\}$ for the scalar product $(.,.)_{sl}$ by the expression
\be
WZW_{sl(m,\C)}=\frac{1}{6}\sum_{i,j,k=1}^{m}(\hat{\tilde{e}}_{i},[\hat{\tilde{e}}_{j},\hat{\tilde{e}}_{k}])_{sl}\n^{R,L}_{\tilde{e}_{i}}\w \n^{R,L}_{\tilde{e}_{j}}\w\n^{R,L}_{\tilde{e}_{k}}.
\ee

\no  Now,  we develop the Wess-Zumino-Witten form of $sl(m\ve n,\C)$ written in the basis $v_{i}$ and dual basis $\hat{v}_{j}$ for $(.,.)_{sl}$ (cf annexe $sl(m\ve n,\C)$) i.e.
\begin{eqnarray}
WZW_{sl(m\ve n,\C)}&=&\sum_{i,j,k=1}^{m+n}(\hat{v}_{i},[\hat{v}_{j},\hat{v}_{k}])_{sl}\n^{R,L}_{v_{i}}\o \n^{R,L}_{v_{j}}\o\n^{R,L}_{v_{k}}\nonumber \\
&=& \frac{1}{6}\sum_{i,j,k=1}^{m+n}(\hat{v}_{i},[\hat{v}_{j},\hat{v}_{k}])_{sl}\n^{R,L}_{v_{i}}\w \n^{R,L}_{v_{j}}\w \n^{R,L}_{v_{k}}\nonumber\\
&=& WZW_{sl(n,\C)}+WZW_{sl(m,\C)}\nonumber\\
&+&(2i)^{3}(\sum_{\g>0}\n^{R,L}_{V_{-\g}}\w\n^{R,L}_{[V_{\g},V_{-\g}]}\w \n_{V_{\g}}\nonumber\\ 
&+&\sum_{\b,\g>0,\b\neq \g,\a}([V_{\b},V_{-\g}],E_{-\a})_{sl}\n^{R,L}_{V_{-\b}}\w\n^{R,L}_{E_{\a}}\w\n^{R,L}_{V_{\g}}\nonumber\\
&+&\sum_{\b,\g>0,\b\neq \g,\a}([V_{\b},V_{-\g}],\tilde{E}_{-\a})_{sl}\n^{R,L}_{V_{-\b}}\w\n^{R,L}_{\tilde{E}_{\a}}\w\n^{R,L}_{V_{\g}}).\nonumber\\
\label{peq}
\end{eqnarray}
With the use of $E_{ts}$, because we have $E_{\a}=E_{st}$, $\tilde{E}_{\b}=iE_{kl}$ and $V_{\g}=E_{ab}$, and the following decomposition
\be
[V_{\g},V_{-\g}]=\sum_{s=1}^{m+n} ([V_{\g},V_{-\g}],\hat{E}_{ss})_{sl}E_{ss}
\ee
the expression $(\ref{peq})$ becomes
\begin{eqnarray}
WZW_{sl(m\ve n,\C)}&=& WZW_{sl(n,\C)}+WZW_{sl(m,\C)}\nonumber \\
&+&(2i)^{3}(\sum_{\ve s \ve=\ve t \ve=0,\ve p \ve =\ve a \ve=1,\ve b \ve =\ve q \ve=0}([E_{ba},E_{pq}],E_{ts})_{sl}\n^{R,L}_{E_{ab}}\w\n^{R,L}_{E_{st}}\w\n^{R,L}_{E_{qp}}\nonumber \\
&-&\sum_{\ve s \ve=\ve t\ve=1,\ve p \ve =\ve a \ve=1,\ve b \ve =\ve q\ve=0 }([E_{ba},E_{pq}],E_{ts})_{sl}\n^{R,L}_{E_{ab}}\w\n^{R,L}_{E_{st}}\w\n^{R,L}_{E_{qp}})\nonumber\\
\label{swzw1}
\end{eqnarray}
On the other hand, we develop
\be
\sum_{s,t,p,q,a,b=1}^{m+n}(\hat{E}_{st},[\hat{E}_{qp},\hat{E}_{ab}])_{sl}\n^{R,L}_{E_{st}}\o \n^{R,L}_{E_{qp}}\o \n^{R,L}_{E_{ab}}
\ee
for this we decompose this sum following the parity of the couples
$(t,s)$, $(p,q)$ and $(a,b)$, and we remark that the term $(\hat{E}_{st},[\hat{E}_{qp},\hat{E}_{ab}])_{sl}$ is often zero \footnote[10]{For instance when $ \hat{E}_{st} $ and $[\hat{E}_{qp},\hat{E}_{ab}] $ have not the same parity.}, then after some algebraic operations we have
\begin{eqnarray}
&&\sum_{s,t,p,q,a,b=1}^{m+n}(\hat{E}_{st},[\hat{E}_{qp},\hat{E}_{ab}])_{sl}\n^{R,L}_{E_{st}}\o \n^{R,L}_{E_{qp}}\o \n^{R,L}_{E_{ab}}\nonumber\\
&=&\frac{1}{6}\sum_{s,t,p,q,a,b=1}^{m}(\hat{E}_{ts},[\hat{E}_{qp},\hat{E}_{ba}])_{sl} \n^{R,L}_{E_{st}}\w \n^{R,L}_{E_{pq}} \w \n^{R,L}_{E_{ab}}\nonumber \\
&+&\frac{1}{6}\sum_{s,t,p,q,a,b=m+1}^{m+n}(\hat{E}_{ts},[\hat{E}_{qp},\hat{E}_{ba}])_{sl} \n^{R,L}_{E_{st}}\w \n^{R,L}_{E_{pq}} \w \n^{R,L}_{E_{ab}}\nonumber \\
&+&(2i)^{3}(\sum_{\ve s \ve=\ve t \ve=0,\ve p \ve =\ve a \ve=1,\ve b \ve =\ve q \ve=0}([E_{ba},E_{pq}],E_{ts})_{sl}\n^{R,L}_{E_{ab}}\w\n^{R,L}_{E_{st}}\w\n^{R,L}_{E_{qp}}\nonumber \\
&-&\sum_{\ve s \ve=\ve t\ve=1,\ve p \ve =\ve a \ve=1,\ve b \ve =\ve q\ve=0 }([E_{ba},E_{pq}],E_{ts})_{sl}\n^{R,L}_{E_{ab}}\w\n^{R,L}_{E_{st}}\w\n^{R,L}_{E_{qp}})\nonumber \\
\label{swzw2}
\end{eqnarray}
Thus, from the equations $(\ref{wzw1})$, $(\ref{wzw2})$, $(\ref{swzw1})$ and $(\ref{swzw2})$ we deduce 
\ben
WZW_{sl(m\ve n,\C)}=\sum_{s,t,p,q,a,b=1}^{m+n}(\hat{E}_{st},[\hat{E}_{qp},\hat{E}_{ab}])_{sl}\n^{R,L}_{E_{st}}\o \n^{R,L}_{E_{qp}}\o \n^{R,L}_{E_{ab}}
\label{swzw3}.
\een
In other words, we have prove the equality $(\ref{egalite})$.


\subsection{Annexe 3: Proof of $C^R=C^L$.}

\no First, we recall the definition of $C^{R,L}$
\be
C^{R,L}=\sum_{a=1}^{m+n}\nabla^{R,L}_{T_{a}}\otimes \nabla^{R,L}_{t_{a}}+(-1)^{\ve a \ve}\nabla^{R,L}_{t_{a}} \otimes \nabla^{R,L}_{T_{a}}.
\ee
From the expression of the basis $T_{a},t_{a}$ (cf. $(\ref{base1})$ and $(\ref{base2})$), we deduce  
\ben
C^{R,L}=C^{R,L}_{u}-C^{R,L}_{w}
\label{eqq1}
\een
with 
\begin{eqnarray}
C^{R,L}_{u,w}&=&2i(\sum_{\alpha>0}(\nabla^{R,L}_{E_{\alpha}}\otimes \nabla^{R,L}_{E_{-\alpha}}+\nabla^{R,L}_{E_{-\alpha}}\otimes \nabla^{R,L}_{E_{\alpha}})+
\sum_{\mu=1}^{m-1}\nabla^{R,L}_{H_{\mu}}\otimes \nabla^{R,L}_{H_{\mu}})\nonumber \\
& &+\sum_{\beta>0}(\nabla^{L,R}_{\tilde{E}_{\beta}}\otimes \nabla^{R,L}_{\tilde{E}_{-\beta}}+\nabla^{R,L}_{\tilde{E}_{-\beta}}\otimes \nabla^{R,L}_{\tilde{E}_{\beta}})+
\sum_{\nu=1}^{n-1}\nabla^{R,L}_{\tilde{H}_{\nu}}\otimes \nabla^{R,L}_{\tilde{H}_{\nu}} \nonumber \\
& &+\sum_{\delta>0}(\nabla^{R,L}_{V_{\delta}} \otimes \nabla^{R,L}_{V_{-\delta}}- 
\nabla^{R,L}_{V_{-\delta}} \otimes \nabla^{R,L}_{V_{\delta}})+ \nabla^{R,L}_{H_{0}} \otimes \nabla^{R,L}_{H_{0}}).
\label{eqq2}
\end{eqnarray}

\no The index $u\;(w)$ means that $C^{R,L}_{u\;(w)}$ acts only on the generators $u(w)$. Moreover, it is not difficult to prove the following equalities
\be
\sum_{\mu=1}^{m-1}H_{\mu} \otimes H_{\mu}=\sum_{k=1}^{m} E_{kk}\otimes E_{kk}-\frac{1}{m}1_{m}\otimes 1_{m}
\ee 
\be
\sum_{\nu=1}^{n-1}\tilde{H}_{\nu} \otimes \tilde{H}_{\nu}=-\sum_{k=m+1}^{m+n} E_{kk}\otimes E_{kk}+\frac{1}{n}1_{n}\otimes 1_{n}
\ee
\be
H_{0}\o H_{0}=\frac{1}{n-m}1_{m+n}\o 1_{m+n}
\ee
where $1_{n}=diag(1_{n},0)$, $1_{m}=diag(0,1_{m})$ and $1_{m+n}=diag(1_{n},1_{m})$. The previous three equations imply
\ben
\sum_{\mu=1}^{m-1}H_{\mu} \otimes H_{\mu}+\sum_{\nu=1}^{n-1}\tilde{H}_{\nu} \otimes \tilde{H}_{\nu} + H_{0} \otimes H_{0}= \sum_{k=1}^{m} E_{kk}\otimes E_{kk}-\sum_{l=m+1}^{m+n}E_{kk}\otimes E_{kk}+\frac{1}{n-m}1 \otimes 1.
\label{eqq3}
\een
Thus, from the equations  $(\ref{eqq2})$ and $(\ref{eqq3})$ we deduce 
\be
C^{R,L}_{u,w}=2i(\sum_{t,s=1}^{m+n}(-1)^{\ve s \ve} \nabla^{R,L}_{E_{st}} \otimes 
\nabla^{R,L}_{E_{ts}}+\frac{1}{n-m}\nabla^{R,L}_{1} \otimes 
\nabla^{R,L}_{1}).
\ee
Then, it is easy to prove respectively the equalities $C^{L}_{u}=C^{R}_{u}$ and $C^{L}_{w}=C^{R}_{w}$ on the generators $u$ and $w$. Therefore, we have $C^{L}_{u}=C^{R}_{u}$ and $C^{L}_{w}=C^{R}_{w}$. Finaly, the equation $(\ref{eqq1})$ and the two previous imply 
\be
C^{R}=C^{L}.
\ee

\section*{Acknowledgements}

It is a pleasure to thank my advisor C.Klimcik for suggesting this 
problem to me, for his invaluable mathematical help and inspiring ideas.



\begin{thebibliography}{00}

\bibitem{A}
Alekseev, A. Yu. \textit{On Poisson actions of compact Lie groups on symplectic manifolds.}  J. Differential Geom. 45 (1997), no. 2, 241--256,

\bibitem{Andru}
Andruskiewitsch, N. \textit{Lie superbialgebras and Poisson-Lie supergroups}.  Abh. Math. Sem. Univ. Hamburg 63 (1993), 147--163,

\bibitem{Ratiu}
Flaschka, H.; Ratiu, T.A. \textit{Convexity theorem for Poisson actions of compact Lie groups.}  Ann. Sci. ƒcole Norm. Sup. (4) 29 (1996), no. 6, 787--809,

\bibitem{GKP}
Grosse, H.; Klim\v{c}\'{\i}k, C.; Pre\v{s}najder, P. \textit{Field theory on a supersymmetric lattice.}  Comm. Math. Phys. 185 (1997), no. 1, 155--175,

\bibitem{Klimcik}
Klimcik, C.; \textit{Quasitriangular WZW model.}  Rev. Math. Phys. 16 (2004), no. 6, 679--808,

\bibitem{LS}
Levendorski\v{i}, S.; Soibelman, Y. \textit{Algebras of functions on compact quantum groups, Schubert cells and quantum tori.}  Comm. Math. Phys. 139 (1991), no. 1, 141--170,

\bibitem{LW}
Lu, J.; Weinstein, A. \textit{Poisson Lie groups, dressing transformations, and Bruhat decompositions.}  J. Differential Geom. 31 (1990), no. 2, 501--526

\bibitem{SNR1} Nahm, W.; Rittenberg, V.; Scheunert, M. \textit{Irreducible representations of the ${\rm osp}(2,1)$ and ${\rm spl}(2,1)$ graded Lie algebras.}  J. Mathematical Phys. 18 (1977), no. 1, 155--162.

\bibitem{SNR2}
Nahm, W.; Rittenberg, V.; Scheunert, M. \textit{Graded Lie algebras: Generalization of Hermitian representations.}  J. Mathematical Phys. 18 (1977), no. 1, 146--154.




\bibitem{Pelleg2}
Pellegrini, F. \textit{Iwasawa decomposition of the Lie supergroup $SL(n,m,\C)$}.  Journal of Algebra and its Applications (accepted for publication),

\bibitem{Pelleg3}
Pellegrini, F. \textit{Les super groupes alg\`ebriques r\'eels et les structures de super Poisson-Lie.} thesis of the author, to appear on arXiv,

\bibitem{Perrin}
Perrin, D. \textit{G\'eom\'etrie alg\`ebrique. Une introduction.} Savoirs Actuels. InterEditions, Paris; CNRS ƒditions, Paris, 1995,

\bibitem{STS}
Semenov-Tian-Shansky, M.A. \textit{Dressing transformations and Poisson group actions.}  Publ. Res. Inst. Math. Sci. 21 (1985), no. 6, 1237--1260, 

\bibitem{Serganova}
Serganova, V.V. \textit{Classification of real simple Lie 
superalgebras and symmetric superspaces}, Functional Analysis {\bf 17 n3} 
(July-September 1983) 46--54,

\bibitem{Waterhouse}
Waterhouse, W. C.  \textit{Introduction to affine group schemes.} Graduate Texts in Mathematics, 66. Springer-Verlag, New York-Berlin, 1979,

\end{thebibliography}
\end{document}